\documentclass{amsart}

\usepackage{graphicx}
\usepackage{amsthm}
\usepackage{amsmath}
\usepackage{amsfonts}
\usepackage{amssymb}

\newtheorem{theorem}{Theorem}

\newtheorem{conjecture}{Conjecture}
\newtheorem*{corollary}{Corollary}
\newtheorem*{definition}{Definition}
\newtheorem*{example}{Example}
\newtheorem*{lemma}{Lemma}
\newtheorem*{notation}{Notation}
\newtheorem*{proposition}{Proposition}
\newtheorem*{remark}{Remark}

\begin{document}
\title[Characteristically nilpotent Lie algebras]{Characteristically nilpotent Lie algebras : a survey}
\author[J.M.Ancochea]{Jos\'e Mar\'{\i}a Ancochea}
\address{Departamento de Geometr\'{i}a y Topolog\'{i}a\\
Facultad CC. Matem\'{a}ticas Univ. Complutense\\
28040 Madrid (Spain)}
\email{Jose\_Ancochea@mat.ucm.es}

\author[O.R.Campoamor]{Rutwig Campoamor}
\address{Departamento de Geometr\'{i}a y Topolog\'{i}a\\
Facultad CC. Matem\'{a}ticas Univ. Complutense\\
28040 Madrid (Spain)}
\email{rutwig@nfssrv.mat.ucm.es}
\thanks{Research partially supported by the D.G.I.C.Y.T project PB98-0758}

\maketitle

\section{Introduction}

The theory of characteristically nilpotent Lie algebras constitutes an
independent research object since 1955. Until then, most studies about Lie
algebras were oriented to the classical aspects of the theory, such as
semisimple and reductive Lie algebras [92]. Though there exists a precedent in
the theory of nilpotent Lie algebras; the Ph.D. thesis of K. Umlauf [95] in
1891, their structure was practically unknown and only classical results like
Engel's theorem were known. From 1939 on, when Lie theorists were seeking from
adequate presentations of the semisimple Lie algebras in terms of generators
and relations, N.\ Jacobson proved that the exceptional complex simple Lie
algebra $G_{2}$ of dimension 14 could be presented as the algebra of
derivations of the Cayley algebra [47]. This result increased the interest in
analyzing the derivations of an arbitrary Lie algebra. However, it was not
until the fifties when the first determining results about derivations of nilpotent
Lie algebras were obtained. It was proven that any nilpotent Lie algebra
has an outer derivation, i.e., there exists at least one derivation which is
not the adjoint operator for a vector of the algebra. Two years earlier, E. V.
Schenkman [85] had published his derivation tower theorem for centerless Lie
algebras, which described in a nice manner the derivation algebras. This
theory was not applicable to the nilpotent algebras, as the adjoint
representation is not faithful. This fact led to the assumption that the
structure of derivations for nilpotent Lie algebras is much more difficult
than for classical algebras. Again, Jacobson proved in 1955 that any Lie
algebra over a field of characteristic zero which has nondegenerate
derivations is nilpotent. In the same paper [48] he asked for the converse.
This result is assumed to be the origin of the theory of characteristically
nilpotent Lie algebras. Dixmier and Lister [29] gave a
negative answer to the converse of Jacobson's theorem. They defined a
generalization of the central descending sequence and called the algebras
satisfying the nullity of a power characteristically nilpotent. The example of
Dixmier and Lister constituted the milestone for a new class of Lie algebras
which seem, in appearance, to be scarce. The first paper about the structure
of characteristically nilpotent Lie algebras, short CNLA, is due to Leger and
T\^{o}g\^{o} in 1959. They proved the equivalence of the sequence condition of Dixmier
and Lister and the nilpotence of the Lie algebra of derivations. Although this
paper does not give any additional example of such an algebra, it reduces the
search to the class of nilpotent Lie algebras. On the other side, the deduced
properties of a CNLA excluded the 2-step nilpotent or metabelian Lie algebras.
The last author, S. T\^{o}g\^{o}, published in 1961 an excellent work which
contained much of the information known about derivation algebras of Lie
algebras. Among others, he introduced special classes of algebras which were shown to be
non CNLA [93]. The importance of CNLAs within the variety of nilpotent Lie
algebra laws was soon recognized by the author, and he also formulated an
interesting question which is nowadays not satisfactorily solved : the problem
of T\^{o}g\^{o}. He asked for the existence of CNLA of derivations, this is,
algebras for which both the derivations and the derivations of these are
nilpotent. Very little is known about the general structure of such Lie
algebras, though its existence has been verified by various authors [8]. The
deformations theory for algebraic structures of M. Gerstenhaber in 1964 [38],
originally developed to study the rigidity of algebraic structures, has become
since then a powerful tool to determine the nilpotence of derivations.
\newline M. Vergne  [96] applied in 1966 the cohomology theory of Lie algebras 
[60] to the study of the variety of nilpotent Lie algebras, obtaining
in particular interesting results about its irreducible components. In
particular, she showed the existence of only two naturally graded filiform Lie
algebras, $L_{n}$ and $Q_{n}$, the second existing only in even dimension. In
particul\~{n}ar, the first algebras has been a central research object for the
last thirty years. Studying its deformations, lots of families of CNLA have
been constructed [52]. In 1970 J. L. Dyer gave a nine dimensional example of
CNLA [33], which was interesting in its own as it had an unipotent
automorphism group. This property is not satisfied by the original example of
Dixmier and Lister, and showed than even CNLA can have quite different
behaviours. By that time, it was perfectly known that such algebras could
exist only from dimension 7 on, as a consequence of the classification in 1958
of the six dimensional algebras [75]. In 1972 G. Favre discovered the lowest
dimensional CNLA known until then [35], which additionally was of the same
nature as Dyer's example. At the same time, G. D. Leger and E. Luks
investigated the metabelian Lie algebras and proved several results about
their rank, and establishing that rank one algebras were given if the
existence of a characteristic ideal containing the derived subalgebra is
assured. These results can be interpreted as a constructive proof that the
original example of 1957 is the known CNLA with lowest characteristic
sequence. The last author applied in 1976 computational methods to prove the
existence of CNLA in any dimension greater or equal to seven. Four years
later, S. Yamaguchi constructed families of CNLA in arbitrary dimension,
constructions that have been completed and generalized in later years [52].
The topological study of the variety of Lie algebra laws led R. Carles [18]
to study the topological properties of CNLA. Among other results he states
that the set of CNLA is constructible for $n\geq7$. For the particular
dimension 7, he also proves that CNLA do not form an open set. Recently [8]
this result has been generalized to any dimension. Another interesting
approach to the CNLA has been deformation theory applied to the Borel
subalgebras of complex simple Lie algebras, like done by Y. B. Khakimdjanov in
1988 to prove that almost all deformations of the nilradical of Borel
subalgebras of complex simple Lie algebras are characteristically nilpotent.
This has shown that these algebras are in fact in abundance within the variety
of nilpotent laws. M. Goze and the last cited author [40] proved, in 1994,
that for any dimension $n\geq9$ an irreducible component of the filiform
variety $\mathcal{F}_{n}$ contains an open set consisting of CNLA.\newline Filiform Lie algebras, specifically the model filiform Lie algebra $L_{n}$, has been also the fundamental source for constructing families of CNLAs. In particular, its cohomology has been calculated, which has allowed to describe its deformations in a precise manner and characterize those deformations which are characteristically nilpotent [54]. Recently, we have turned our interest to nilpotent Lie algebras which structurally "look like $Q_{n}$". As known, this algebra cannot appear in odd dimension. This is a consequence of the so called centralizer property [22], which codifies information about the structure of the commutator subalgebra and the ideals of the central descending sequence. Now the centralizer property can be generalized to any naturally graded nilpotent Lie algebra, and defines a class of algebras which can be interpreted as those which are the "easiest nilpotent Lie algebras to deform for obtaining CNLAs". The key to this is extension theory combined with deformation theory. \newline This approach also leads to certain questions about the rigidity of a nilpotent Lie algebra. In 1970 Vergne postulated the nonexistence of nilpotent Lie algebras that are rigid in the variety $\frak{L}^{n}$ for $n\neq 1$. In his study about the structure of rigid Lie algebras [18], Carles established that if a nilpotent Lie algebra is rigid, then it necessarily must be a CNLA. The strongness of this result seems to confirm the validity of the conjeture, although there is no known procedure to prove it.\newline Finally, we review some results about affine structures over Lie algebras. This kind of structures are of great importance not only for purposes of cohomology theory [15], but also for representation theory of nilpotent Lie algebras. The interesting point is that CNLA can admit an affine structure, such as it was proven for the example of Dixmier and Lister by Scheunemann [87] in 1974. Although practically nothing is known about CNLA with affine structures, the cohomological method developed by Burde in [15] could be an important source for studying these algebras.   

\section{Generalities}

In this section we resume the elementary facts about Lie algebras that will be
used thorughout the paper. Although it is often unnecessary to specify the
base field, we will asume here that all Lie algebras are complex.

\begin{definition}
Let $\frak{g}$ be a finite dimensional vectorial space over $\mathbb{C}$.
A Lie algebra law over $\mathbb{C}^{n}$ is a bilinear alternated mapping $\mu\in Hom\left(  \mathbb{C}^{n}\times\mathbb{C}^{n},\mathbb{C}^{n}\right)$ which satisfies the conditions

\begin{enumerate}
\item $\mu\left(  X,X\right)  =0,\;\forall\;X\in\mathbb{C}^{n}$

\item $\mu\left(  X,\mu\left(  Y,Z\right)  \right)  +\mu\left(  Z,\mu\left(
X,Y\right)  \right)  +\mu\left(  Y,\mu\left(  Z,X\right)  \right)
=0,\;\forall\;X,Y,Z\in\mathbb{C}^{n}$, \ \newline  ( Jacobi identity )
\end{enumerate}
If  $\mu$ is a Lie algebra law, the pair  $\frak{g}= \left(\mathbb{C}^{n}, \mu \right)$ is  called Lie algebra. From now on we identify the Lie algebra with its law $\mu$. 
\end{definition}

\begin{remark}
We say that $\mu$ is the law of $\frak{g}$, and where necessary we use the bracket notation to describe the law :
\[
\left[  X,Y\right]  =\mu\left(  X,Y\right)  ,\;\forall\;X,Y\in\frak{g}%
\]
The nondefined brackets are zero or obtained by symmetry. 
\end{remark}

\begin{definition}
Given an ideal $\frak{I}$  of $\frak{g}$, we call centralizer of $\frak{I}$
in  $\frak{g}$ to the subalgebra 
\[
C_{\frak{g}}\frak{I}=\left\{  X\in\frak{g}\;|\;\mu\left(  X,\frak{I}\right)
=0\right\}
\]
\end{definition}

\bigskip To any Lie algebra we can associate the two following sequences :

\begin{align*}
D^{0}\frak{g}  &  =\frak{g}\supset D^{1}\frak{g}=\left[  \frak{g}%
,\frak{g}\right]  \supset..\supset D^{k}\frak{g}=\left[  D^{k-1}\frak{g},D^{k-1}\frak{g}\right]
  \supset...\\
C^{0}\frak{g}  &  =\frak{g}\supset C^{1}\frak{g}=D^{1}\frak{g}\supset
C^{2}\frak{g}=\left[  C^{1}\frak{g},\frak{g}\right]  \supset...\supset
C^{k}\frak{g}=\left[  C^{k-1}\frak{g},\frak{g}\right]  \supset...
\end{align*}
called respectively derived and descending central sequences of $\frak{g}$.

\begin{definition}
Let $\frak{g}$ be a Lie algebra. We say that 

\begin{enumerate}
\item $\frak{g}$ is solvable if there exists an integer $k\geq1$ such that
$D^{k}\frak{g}=\left\{  0\right\}  $.

\item $\frak{g}$ is nilpotent if there exists an integer ( called nilindex $n\left( \frak{g}\right)$ of $\frak{g}$) $k\geq1$ such that  $C^{k}\frak{g}=\left\{
0\right\}  .$
\end{enumerate}
\end{definition}

\begin{definition}
An $n$-dimensional nilpotent Lie algebra is called filiform if 
\[
\dim C^{k}\frak{g}=n-k-1,\;1\leq k\leq n-1
\]
\end{definition}

\begin{remark}
Calling $p_{i}=dim\left(\frac{C^{i-1}\frak{g}}{C^{i}\frak{g}}\right)$ for $1\leq p_{i}\leq n\left( \frak{g}\right)$, the type of the nilpotent Lie algebra is the sequence $\left\{p_{1},..,p_{r}\right\}$. Then a filiform algebra corresponds to those of type $\left\{2,1,..,1\right\}$ [97].
\end{remark}

We recall the laws for the $\left(n+1\right)$-dimensional filiform Lie algebras $L_{n}$ and $Q_{n}$, which are basically the only filiform Lie algebras we have to deal with here :

\begin{enumerate}
\item $L_{n}$ $\left(  n\geq3\right)  :$
\[
\lbrack X_{1},X_{i}]=X_{i+1},\;2\leq i\leq n
\]
over the basis $\left\{X_{1},..,X_{n+1}\right\}$.

\item $Q_{2m-1}$ $\left(  m\geq3\right)  :$%
\begin{align*}
\lbrack X_{1},X_{i}]  & =X_{i+1},\;2\leq i\leq2m-1\\
\lbrack X_{j},X_{2m+1-j}]  & =\left(  -1\right)  ^{j}X_{2m},\;1\leq j\leq m
\end{align*}
over the basis $\left\{X_{1},..,X_{2m}\right\}$.
\end{enumerate}

\begin{definition}
A Lie algebra $\frak{g}$ \ is graded over  $\mathbb{Z}$ if
it admits a decomposition 
\[
\frak{g}=\bigoplus_{k\in\mathbb{Z}}\frak{g}_{k}%
\]
where the $\frak{g}_{k}$ are $\mathbb{C}$-subspaces of
$\frak{g}$ which satisfy  $\left[  \frak{g}_{r},\frak{g}_{s}\right]
\subset\frak{g}_{r+s}$ , $\;r,s\in\mathbb{Z}$.
\end{definition}

Observe that any graduation defines a sequence
\[
S_{k}=F_{k}\left(  \frak{g}\right)  =\bigoplus_{t\geq k}\frak{g}_{t}%
\]
with the properties

\begin{enumerate}
\item $\frak{g}=%
{\displaystyle\bigsqcup}
\,S_{k}$

\item $\left[  S_{i},S_{j}\right]  \subset S_{i+j}\;\forall i,j$

\item $S_{i}\subset S_{j}$ \ si $i>j$
\end{enumerate}

\begin{definition}
A family $\left\{  S_{i}\right\}  $ of subspaces of 
$\frak{g}$ define a filtration ( descending ) over $\frak{g}$ if
it satisfies properties $1), 2), 3)$. The algebra is called filtered.
\end{definition}

The construction can be reversed, i.e., any filtration defines a graduation by taking  $\frak{g}
_{k}=\frac{S_{k-1}}{S_{k}}$ \ for $k\geq1$. The graduation is called associated to the filtration  $\left\{  S_{i}\right\}  $ and it defines a Lie algebra. 

\begin{definition}
A nilpotent Lie algebra is called naturally graded if 
$\frak{g}\simeq\frak{gr}\left(  \frak{g}\right)$, where $\frak{gr}\left(\frak{g}\right)$ is the graduation associated to the filtration induced in $\frak{g}$ by the central descending sequence.
\end{definition}

It follows immediately that both $L_{n}$ and $Q_{n}$ are naturally graded. They are in fact the only filiform Lie algebras having this property [97].

\begin{definition}
A derivation $f$ of a Lie algebra $\frak{g}$ is a linear mapping
\[
f:\frak{g}\longrightarrow\frak{g}%
\]
satisfying
\[
\left[  f\left(  X\right)  ,Y\right]  +\left[  X,f\left(  Y\right)  \right]
-f\left[  X,Y\right]  =0,\;\;\forall\left(  X,Y\right)  \in\frak{g}^{2}%
\]
\end{definition}

\vspace{5mm} \noindent We denote by $Der\frak{g}$ the set of derivations of
$\frak{g}$.\ It is a Lie subalgebra of $End\frak{g}$.
\begin{proposition}
For all $X$ in $\frak{g}$, the endomorphism $adX$ is a derivation of
$\frak{g}$.
\end{proposition}

\begin{definition}
\noindent The derivations $f$ of $\frak{g}$ which are of type $f=adX$ for
$X\in\frak{g}$ are called inner derivations.
\end{definition}

\subsection{Cohomology of Lie algebras}

There exists a general study of the cohomology of Lie algebra by considering
the cohomology with values on a $\frak{g}$-module.\ See for example references
[60].\newline Let $\frak{g}$ be a Lie algebra. A $p$-dimensional
cochain of $\frak{g}$ (with values in $\frak{g}$) is a $p$-linear alternating
mapping of $\frak{g}^{p}$ in $\frak{g}$ $\left(  p\in\mathbb{N}^{\ast}\right)
$. A $0$-cochain is a constant function from $\frak{g}$ to $\frak{g}%
$.\newline We denote by $C^{p}\left(  \frak{g},\frak{g}\right)  $ as the space
of the $p$-cochains and
\[
C^{\ast}\left(  \frak{g},\frak{g}\right)  =\oplus_{p\geq0}C^{p}\left(
\frak{g},\frak{g}\right)  .
\]
We can provide $C^{p}\left(  \frak{g},\frak{g}\right)  $ of a $\frak{g}%
$-module structure by putting
\[
\left(  X\Phi\right)  \left(  X_{1},...,X_{p}\right)  =\left[  X,\Phi\left(
X_{1},...,X_{p}\right)  \right]  -\sum_{1\leq i\leq p}\Phi\left(
X_{1},...,\left[  X,X_{i}\right]  ,...,X_{p}\right)
\]
for all $\ X_{1},...,X_{p}\ \in\frak{g}.$

On the space $C^{\ast}\left(  \frak{g},\frak{g}\right)  $ we define the
endomorphism
\[
\delta:C^{\ast}\left(  \frak{g},\frak{g}\right)  \longrightarrow C^{\ast
}\left(  \frak{g},\frak{g}\right)
\]%
\[
\Phi\longrightarrow\delta\Phi
\]
by putting
\[
\delta\Phi\left(  X\right)  =X.\Phi\ \ \ \mbox{if}\ \ \ \Phi\in C^{0}\left(
\frak{g},\frak{g}\right)
\]%
\[
\delta\Phi\left(  X_{1},...,X_{p+1}\right)  =\sum_{1\leq s\leq p+1}\left(
-1\right)  ^{s+1}\left(  X_{s}.\Phi\right)  \left(  X_{1},...,\overset
{}{\widehat{X}}_{s},...,X_{p+1}\right)  +
\]%
\[
+\sum_{1\leq s\leq t\leq p+1}\left(  -1\right)  ^{s+t}\Phi\left(  \left[
X_{s},X_{t}\right]  ,X_{1},...,\overset{}{\widehat{X}_{s}},...,\overset
{}{\widehat{X}}_{t,...},X_{p+1}\right)
\]
if $\Phi\in C^{p}\left(  \frak{g},\frak{g}\right)  $, $p\geq1$.

\noindent By this definition, $\delta\left(  C^{p}\left(  \frak{g}%
,\frak{g}\right)  \right)  \subset C^{p+1}\left(  \frak{g},\frak{g}\right)  $
and we can verify that
\[
\delta\circ\delta=0.
\]

We denote by
\[
\left\{
\begin{array}
[c]{cc}%
Z^{p}\left(  \frak{g},\frak{g}\right)  =Kerd\left|  _{C^{p}\left(
\frak{g},\frak{g}\right)  }\right.   & p\geq1\\
B^{p}\left(  \frak{g},\frak{g}\right)  =Imd\left|  _{C^{p}\left(
\frak{g},\frak{g}\right)  }\right.   & p\geq1
\end{array}
\right.
\]
and $H^{p}\left(  \frak{g},\frak{g}\right)  =Z^{p}\left(  \frak{g}%
,\frak{g}\right)  \;|\;B^{p}\left(  \frak{g},\frak{g}\right)  ,$ $p\geq1$.

\noindent This quotient space is called the cohomology space of $\frak{g}$ of
degree $p$ (with values in $\frak{g}$). For $p=0$, then we put $B^{0}\left(
\frak{g},\frak{g}\right)  =\left\{  0\right\}  $ and $H^{0}\left(
\frak{g},\frak{g}\right)  =Z^{0}\left(  \frak{g},\frak{g}\right)  $. This last
space can be identified to the space of all $\frak{g}$-invariant elements that
is
\[
\left\{  X\in\frak{g}\ \mbox{such that}\ adY\left(  X\right)  =0\ \ \forall
Y\in\frak{g}\right\}  .
\]
Then $Z^{0}\left(  \frak{g},\frak{g}\right)  =Z\left(  \frak{g}\right)  $ (the
center of $\frak{g}$).

\subsubsection{The space $H^{1}\left(  \frak{g},\frak{g}\right)  $}

We have
\[
Z^{1}\left(  \frak{g},\frak{g}\right)  =\left\{  f:\frak{g}\longrightarrow
\frak{g\;}|\;\delta f=0\right\}  .
\]
But $\delta f\left(  X,Y\right)  =\left[  f\left(  X\right)  ,Y\right]
+\left[  X,f\left(  Y\right)  \right]  -f\left[  X,Y\right]  $.\ Then
$Z^{1}\left(  \frak{g},\frak{g}\right)  $ is nothing but the algebra of
derivation of $\frak{g}$:
\[
Z^{1}\left(  \frak{g},\frak{g}\right)  =Der\frak{g}.
\]
It is the same for :
\[
B^{1}\left(  \frak{g},\frak{g}\right)  =\left\{  adX,X\in\frak{g}\right\}  .
\]
Thus the space $H^{1}\left(  \frak{g},\frak{g}\right)  $ can be interpreted as
the set of the outer derivations of the Lie algebra $\frak{g}.$

\bigskip Let $I$ be an ideal of $\frak{g}$.\ We consider the cochains
\[
\varphi:I^{p}\longrightarrow\frak{g}%
\]
on $I$ with values in $\frak{g}$. For these cochains we can also define, by
restriction, the coboundary operator $\delta$. As $I$ is an ideal of
$\frak{g}$, $H^{1}\left(  \frak{g},\frak{g}\right)  $ is a $\frak{g}$-module.
So we can consider the cohomology space $H^{*}\left(  I,\frak{g}\right)  $.

\noindent A $p$-cochain $\varphi$ of $C^{p}\left(  I,\frak{g}\right)  $ is
$\frak{g}$-invariant if it satisfies :
\[
X\varphi\left(  X_{1},...,X_{p}\right)  =\left[  X,\varphi\left(
X_{1},...,X_{p}\right)  \right]  -{\sum_{1\leq i\leq p}}\varphi\left(
X_{1},...,X_{i-1},\left[  X,X_{i}]\right]  ,...,X_{p}\right)  =0
\]
We denote by $C^{\ast}\left(  I,\frak{g}\right)  ^{\frak{g}}$ the set of
cochains on $I$ which are $\frak{g}$-invariant and $H^{\ast}\left(
I,\frak{g}\right)  ^{\frak{g}}$ the correspondent cohomology space. Each
element $\overline{\varphi}$ of $H^{p}\left(  I,\frak{g}\right)  ^{\frak{g}}$
has a representative which is the restriction to $I$ of a cochain $\Psi$ in
$C^{p}\left(  \frak{g},\frak{g}\right)  $ such that $d\psi\in\left(
\frak{g}/I,\frak{g}^{I}\right)  $ where $\frak{g}^{I}=\left\{  X\in
\frak{g}/\left[  X,Y\right]  =0\ \ \forall Y\in I\right\}  $. This element
$d\psi$ does not depend upon the choice of the representative of
$\overline{\varphi}$. Let $t_{p+1}$be the homomorphism so defined :
\[
t_{p+1}:H^{p}\left(  I,\frak{g}\right)  ^{\frak{g}}\longrightarrow
H^{p+1}\left(  \frak{g}/I,\frak{g}^{I}\right)  .
\]
We define an exact sequence :
\[
0\longrightarrow H^{p}\left(  \frak{g}/I,\frak{g}^{I}\right)  \overset{l_{p}%
}{\longrightarrow}H^{p}\left(  \frak{g},\frak{g}\right)  \overset{r_{p}%
}{\longrightarrow}H^{p}\left(  I,\frak{g}\right)  ^{\frak{g}}%
\]%
\[
\overset{r_{p+1}}{\longrightarrow}H^{p+1}\left(  \frak{g}/I,\frak{g}%
^{I}\right)  \overset{}{\longrightarrow}H^{p+1}\left(  \frak{g},\frak{g}%
\right)
\]
where $r_{p}$ is the homomorphism restriction and $l_{p}$ is defined by
looking upon the cochains of $\frak{g}/I$ in $\frak{g}^{I}$ as cochain of
$\frak{g}$ in $\frak{g}$.

\begin{example}
We suppose that $codim\left(  I\right)  =1$.\noindent Then $\dim\frak{g}/I=1$
and $C^{p}\left(  \frak{g}/I,\frak{g}\right)  =0$ for $p\geq2$. Thus
\[
0\longrightarrow0\longrightarrow H^{2}\left(  \frak{g},\frak{g}\right)
\longrightarrow H^{2}\left(  I,\frak{g}\right)  ^{\frak{g}}\longrightarrow0
\]
and we have
\[
H^{2}\left(  \frak{g},\frak{g}\right)  =H^{2}\left(  I,\frak{g}\right)
^{\frak{g}}.
\]
\end{example}

\subsection{The spaces $H^{2}\left(\frak{g},\mathbb{C}\right)$}

Recall that the space $H^{2}\left(  \frak{g},\mathbb{C}^{p}\right)  $
can be interpreted as the space of classes of $p$-dimensional central extensions of the Lie algebra $\frak{g}$. We recall the elementary facts :

Let $\frak{g}$ be an $n$-dimensional nilpotent Lie algebra with law $\mu_{0}$. 
A central extension of $\frak{g}$ by $\mathbb{C}^{p}$ is an exact sequence of Lie algebras
\[
0\longrightarrow\mathbb{C}^{p}\longrightarrow\overset{-}{\frak{g}%
}\longrightarrow\frak{g}\longrightarrow0\text{ }%
\]
such that $\mathbb{C}^{p}\subset Z\left(  \overset{-}{\frak{g}}\right)  .$ Let
$\alpha$ be a cocycle of the De Rham cohomology $Z^{2}\left(  \frak{g},\mathbb{C}%
^{p}\right)  .$ This gives the extension
\[
0\longrightarrow\mathbb{C}^{p}\longrightarrow\mathbb{C}^{p}\oplus\frak{g}%
\longrightarrow\frak{g}\longrightarrow0
\]
with associated law $\mu=\mu_{0}+\alpha$ defined by%
\[
\mu\left(  \left(  a,x\right)  ,\left(  b,y\right)  \right)  =\left(
\alpha\mu_{0}\left(  x,y\right)  ,\mu_{0}(x,y\right)  )
\]
In the following we are only interested in extensions of $\mathbb{C}$ by $\frak{g}$, i.e, extensions of degree one. It is well known that the space of $2$-cocycles
$Z^{2}\left(  \frak{g},\mathbb{C}\right)$ is identified with the space of linear forms
over $\bigwedge^{2}\frak{g}$ which are zero over the subspace $\Omega$ :%
\[
\Omega:=\left\langle \mu_{0}\left(  x,y\right)  \wedge z+\mu_{0}\left(
y,z\right)  \wedge x+\mu_{0}\left(  z,x\right)  \wedge y\right\rangle
_{\mathbb{C}}%
\]
The extension classes are defined modulus the coboundaries 
$B^{2}\left(  \frak{g},\mathbb{C}\right)  .$ This allows to identify the cohomology space
$H^{2}\left(  \frak{g},\mathbb{C}\right)  $ with the dual of the space 
$\frac{Ker\;\lambda}{\Omega},$ where $\lambda\in Hom\left(  \bigwedge
^{2}\frak{g},\frak{g}\right)  $ is defined as
\[
\lambda\left(  x\wedge y\right)  =\mu_{0}\left(  x,y\right)\; x,y\in\frak{g}
\]
In fact we have $H_{2}\left(\frak{g},\mathbb{C}\right)=\frac{Ker\lambda}{\Omega}$ for the $2$-homology space, and as $H^{2}\left(\frak{g},\mathbb{C}\right)=Hom_{\mathbb{C}}\left(H_{2}\left(\frak{g},\mathbb{C}\right),\mathbb{C}\right)$ the assertion follows.

\begin{notation}
Let $\varphi_{ij}\in H^{2}\left(  \frak{g},\mathbb{C}\right)  $ the cocycles
defined by 
\[
\varphi_{ij}\left(  X_{k},X_{l}\right)  =\delta_{ik}\delta_{jl}
\]
\end{notation}

Observe that a cocycle  $\varphi$ can be written as a linear combination of the preceding cocycles. We have  :

\begin{lemma}
$\sum a^{ij}\varphi_{ij}=0$ if and only if  $\sum a^{ij}\left(  X_{i}\wedge
X_{j}\right)\in\Omega$
\end{lemma}

Let $\frak{g}$ be an $n$-dimensional nilpotent Lie algebra. The subspace of central extensions is noted by  $E_{c,1}\left(  \frak{g}\right)$. It has been shown that this space is 
irreducible and constructible. However, for our purpose this space is too general. We only need certain cohomology classes of this space.

\begin{notation}
For $k\geq2$ let %
\[
H_{k}^{2,t}\left(  \frak{g},\mathbb{C}\right)  =\left\{  \varphi_{ij}\in
H^{2}\left(  \frak{g},\mathbb{C}\right)  \;|\;i+j=2t+1+k\right\}  ,\;1\leq t\leq
\left[  \frac{n-3}{2}\right],\
\]

\[
H_{k}^{2,\frac{t}{2}}\left(  \frak{g},\mathbb{C}\right)  =\left\{
\varphi_{ij}\in H^{2}\left(  \frak{g},\mathbb{C}\right)
\;|\;i+j=t+1+k\right\}  ,\;t\in\{1,..,\left[  \frac{n-3}{2}\right]
\},t\equiv1\left(  \operatorname{mod}\,2\right)
\]

\end{notation}

These cocycles are essential to determine the central extensions which are additionally naturally graded. If $\mathbf{E}_{c,1}\left(\frak{g}\right)  $ denotes the central extensions that are naturally graded, we consider the subspaces 
\[
\mathbf{E}_{c,1}^{t,k_{1},..,k_{r}}\left(  \frak{g}\right)  =\left\{  \mu
\in\mathbf{E}_{c,1}\left(  \frak{g}\right)  \;|\;\mu=\mu_{0}+\left(
{\textstyle\sum}
\varphi_{ij}^{k_{i}}\right)  ,\;\varphi_{ij}^{k_{i}}\in H_{k_{i}}^{2,t}\left(
\frak{g},\mathbb{C}\right)  \right\}
\]%

\[
\mathbf{E}_{c,1}^{\frac{t}{2},k_{1},..,k_{r}}\left(  \frak{g}\right)
=\left\{  \mu\in\mathbf{E}_{c,1}\left(  \frak{g}\right)  \;|\;\mu=\mu
_{0}+\left(
{\textstyle\sum}
\varphi_{ij}^{k_{i}}\right)  ,\;\varphi_{ij}^{k_{i}}\in H_{k_{i}}^{2,\frac{t}{2}%
}\left(  \frak{g},\mathbb{C}\right)  \right\}
\]
where $0\leq k_{j}\in\mathbb{Z},$ $j=1,..,r.$\newline
Given a basis $\left\{X_{1},..,X_{n},X_{n+1}\right\}$ of $\mu$ belonging to any of these spaces, the Lie algebra law is defined by : 

\[
\mu\left(  X_{i},X_{j}\right)  =\mu_{0}\left(  X_{i},X_{j}\right)  +\left(
\sum\varphi_{ij}^{k}\right)  X_{n+1},\;1\leq i,j\leq n\; \left(X_{i},X_{j}\right)\in\frak{g}^{2}
\]

\begin{lemma}
As vector spaces, the following identity holds :
\[
\mathbf{E}_{c,1}\left(  \frak{g}\right)  =\sum_{t,k}\mathbf{E}_{c,1}%
^{t,k_{1},..,k_{r}}\left(  \frak{g}\right)  +\mathbf{E}_{c,1}^{\frac{t}%
{2},k_{1},..,k_{r}}\left(  \frak{g}\right)
\]
\end{lemma}

This follows easily. Observe that, though $t$ is bounded by the dimension, $k\geq2$
has no restrictions. However, the sum is finite, for the spaces $\mathbf{E}_{c,1}^{t,k_{1},..,k_{r}}$  are zero
for almost any choice $\left( k_{1},.., k_{r}\right)$.\newline Given the Lie algebra  $\frak{g}=\left(\mathbb{C}^{n},\mu_{0}\right)$, we have the associated graduation  
$\frak{gr}\left(  \frak{g}\right)  =\sum_{i=1}^{n\left(  \frak{g}\right)
}\frak{g}_{i},$ where $\frak{g}_{i}=\frac{C^{i-1}\frak{g}}{C^{i}\frak{g}}$ and
$n\left(  \frak{g}\right)  $ is the nilindex of $\frak{g}.$
Independently
of $\frak{g}$ being naturally graded or not, any
vector $X$ has a fixed position in one of the graduation blocks.
The study of the central extensions which preserve a graduation is reduced to the study of the position of the adjoined vector $X_{n+1}$. Note that in this sense the cocycles
$\varphi_{ij}\in H_{k}^{2,t}\left(  \frak{g},\mathbb{C}\right)  $ codify 
this information.

\subsection{The algebraic variety $\frak{L}^{n}$.}
A $n$-dimensional complex Lie algebra can be seen as a pair $\frak{g}%
=(\mathbb{C}^{n},\mu)$ where $\mu$ is a Lie algebra law on $\mathbb{C}^{n},$
the underlying vector space to $\frak{g}$ is $\mathbb{C}^{n}$ and $\mu$ the
bracket of $\frak{g}$. We will note by $\frak{L}^{n}$ the set of Lie algebra laws on
$\mathbb{C}^{n}$. It is a subset of the vectorial space of alternating
bilinear mappings on $\mathbb{C}^{n}.$

\begin{definition}
Two laws $\mu$ and $\mu^{\prime}\in \frak{L}^{n}$ are said isomorphic, if there is
$f\in Gl(n,\mathbb{C)}$ such that
\[
\mu^{\prime}(X,Y)=f*\mu(X,Y)=f^{-1}(\mu(f(X),f(Y)))
\]
for all $X,Y\in\mathbb{C}^{n}.$
\end{definition}

In this case, the Lie algebras $\frak{g}=(\mathbb{C}^{n},\mu)$ and
$\frak{g}^{\prime}=(\mathbb{C}^{n},\mu^{\prime})$ are isomorphic.

Let $\mathcal{O(\mu)}$ be the set of the laws isomorphic to $\mu.$ It is
called the orbit of $\mu$.

Let us fix a basis $\{e_{1},e_{2},\cdots,e_{n}\}$ of $\mathbb{C}^{n}.$ The
structural constants of $\mu\in \frak{L}^{n}$ are the complex numbers $C_{ij}^{k}$
given by
\[
\mu(e_{i},e_{j})=\sum_{k=1}^{n}C_{ij}^{k}e_{k}.
\]
As the basis is fixed, we can identify the law $\mu$ with its structural
constants. These constants satisfy :
\[
(1)\left\{
\begin{array}
[c]{ll}%
C_{ij}^{k}=-C_{ji}^{k}\ ,\ \ \ 1\leq i<j\leq n\ ,\ \ \ 1\leq k\leq n & \\
\sum_{l=1}^{n}C_{ij}^{l}C_{lk}^{s}+C_{jk}^{l}C_{li}^{s}+C_{ki}^{l}C_{jl}%
^{s}=0\ ,\ \ \ 1\leq i<j<k\leq n\ ,\ \ \ 1\leq s\leq n. &
\end{array}
\right.
\]
Then $\frak{L}^{n}$ appears as an algebraic variety embedded in the linear space of
alternating bilinear mapping on $\mathbb{C}^{n}$, isomorphic to $\mathbb{C}%
^{\frac{n^{3}-n^{\acute{e}}}{2}}.$

Let be $\mu\in \frak{L}^{n}$ and consider the Lie subgroup $G_{\mu}$ of
$Gl(n,\mathbb{C)}$ defined by
\[
G_{\mu}=\{f\in Gl(n,\mathbb{C)}\mid\text{\quad}f*\mu=\mu\}
\]
Its Lie algebra is the Lie algebra of derivations of $\mu.$ Let be
$\mathcal{O}(\mu)$ the orbit of $\mu$ respect the action of $Gl(n,\mathbb{C)}%
$. It is isomorphic to the homogeneous space $Gl(n,\mathbb{C)}$/$G_{\mu}.\,
$Then it is a $\mathcal{C}^{\infty}$ differential manifold of dimension
\[
\dim\mathcal{O}(\mu)=n^{2}-\dim Der(\mu).
\]

It is not difficult to see that the orbit $\mathcal{O}$($\mu$) of $\mu$ is a differentiable
manifold [96] embedded in $\frak{L}^{n}$ defined by
\[
\mathcal{O}(\mu)=\frac{Gl(n,\mathbb{C)}}{G_{\mu}}%
\]
We consider a point $\mu^{\prime}$ close to $\mu$ in $\mathcal{O}$($\mu$).
There is $f\in Gl(n,\mathbb{C)}$ such that $\mu^{\prime}=f*\mu.$ Suppose that
$f$ is close to the identity : $f=Id+\varepsilon g,$ with $g\in gl(n)$ Then
\begin{align*}
\mu^{\prime}(X,Y)  &  =\mu(X,Y)+\varepsilon[-g(\mu(X,Y))+\mu(g(X),Y)+\mu
(X,g(Y))]\\
&  +\varepsilon^{2}[\mu(g(X),g(Y))-g(\mu(g(X),Y)+\mu(X,g(Y))-g\mu(X,Y)].
\end{align*}
Then
\[
\frac{\mu^{\prime}(X,Y)-\mu(X,Y)}{\varepsilon}\rightarrow_{\varepsilon
\rightarrow0}\delta_{\mu}g(X,Y)
\]
Among the possible orbits, those which are open are specially important for the study of the variety, as we will see later.

\begin{definition}
Let $\mu$ be a law such that the orbit $\mathcal{O}$($\mu)$ is open in $\frak{L}^{n}$. Then $\mu$ is called a rigid law.
\end{definition}

\begin{proposition}
The tangent space to the orbit $\mathcal{O}$($\mu)$ at the point $\mu$ is the
space $B^{2}(\mu,\mu)$ of the 2-cocycles of the Chevalley cohomology of $\mu.$
\end{proposition}

Let $\mu$ be in $\frak{L}^{n}$ and consider a bilinear alternating mapping
$\mu^{\prime}=\mu+t\varphi$ where $t$ is a small parameter. Then $\mu^{\prime
}\in \frak{L}^{n}$ for all $t$ if and only if we have :
\[
\left\{
\begin{array}
[c]{c}%
\delta\varphi=0\\
\varphi\in \frak{L}^{n}%
\end{array}
\right.
\]

\begin{proposition}
A straight line $\Delta$ passing throught $\mu$ is a tangent line in $\mu$ to
$\frak{L}^{n}$ if its direction is given by a vector of $Z^{2}(\mu,\mu).$
\end{proposition}

Suppose that $H^{2}(\mu,\mu)=0.$ Then the tangent space to $\mathcal{O}$%
($\mu)$ at the point $\mu$ is the set of the tangent lines to $\frak{L}^{n}$ at the
point $\mu.$ Thus the tangent space to $\frak{L}^{n}$ exists in this point and it is
equal to $B^{2}(\mu,\mu).$ The point $\mu$ is a nonsingular point. We deduce
of this that the inclusion $\mathcal{O}$($\mu)\hookrightarrow \frak{L}^{n}$ is a
local homeomorphism. This property is valid for all points of $\mathcal{O}%
$($\mu)$, then $\mathcal{O}$($\mu)$ is open in $\frak{L}^{n}$ (for the induced metric topology).

\begin{proposition}
Let $\mu\in \frak{L}^{n}$ such that $H^{2}(\mu,\mu)=0.$ If the algebraic variety
$\frak{L}^{n}$ is provided with the metric topology induced by $\mathbb{C}%
^{\frac{n^{3}-n^{2}}{2}}$, then the orbit $\mathcal{O}$($\mu)$ is open in
$\frak{L}^{n}$.
\end{proposition}

This geometrical approach shows the problems undelying to the existence of
singular points in the algebraic variety $\frak{L}^{n}$ [21].

\subsection{Formal deformations}

Let be $\varphi,\psi\in C^{2}(\mathbb{C}^{n},\mathbb{C}^{n})$ two
skew-symmetric bilinear maps on $\mathbb{C}^{n}.$ We define the trilinear
mapping $\varphi\circ\psi$ on $\mathbb{C}^{n}$ by
\[
\varphi\circ\psi(X,Y,Z)=\varphi(\psi(X,Y),Z)+\varphi(\psi(Y,Z),X)+\varphi
(\psi(Z,X),Y)
\]
for all $X,Y,Z\in\mathbb{C}^{n}.$ Using this notation, the Lie bracket is
written $\mu\circ\mu=0.$

Let be $\mu_{0}\in L^{n}$ and $\varphi\in C^{2}(\mathbb{C}^{n},\mathbb{C}%
^{n}).$ Then $\varphi\in Z^{2}(\mu_{0},\mu_{0})$ if and only if
\[
\mu_{0}\circ\varphi+\varphi\circ\mu_{0}=\delta_{\mu_{0}}\varphi=0.
\]

\begin{definition}
A (formal) deformation of a law $\mu_{0}\in L^{n}$ is a formal sequence with
parameter $t$%
\[
\mu_{t}=\mu_{0}+\sum_{t=1}^{\infty}t^{i}\varphi_{i}%
\]
where the $\varphi_{i}$ are skew-symmetric bilinear maps $\mathbb{C}^{n}%
\times\mathbb{C}^{n}\rightarrow\mathbb{C}^{n}$ such that $\mu_{t}$ satisfies
the formal Jacobi identity $\mu_{t}\circ\mu_{t}=0.$
\end{definition}

Let us develop this last equation.
\[
\mu_{t}\circ\mu_{t}=\mu_{0}\circ\mu_{0}+t\delta_{\mu_{0}}\varphi_{1}%
+t^{2}(\varphi_{1}\circ\varphi_{1}+\delta_{\mu_{0}}\varphi_{2})+t^{3}%
(\varphi_{1}\circ\varphi_{2}+\varphi_{2}\circ\varphi_{1}+\delta_{\mu_{0}%
}\varphi_{3})+...
\]
and the formal equation $\mu_{t}\circ\mu_{t}=0$ is equivalent to the infinite
system
\[
(I)\text{ }\left\{
\begin{array}
[c]{l}%
\mu_{0}\circ\mu_{0}=0\\
\delta_{\mu_{0}}\varphi_{1}=0\\
\varphi_{1}\circ\varphi_{1}=-\delta_{\mu_{0}}\varphi_{2}\\
\varphi_{1}\circ\varphi_{2}+\varphi_{2}\circ\varphi_{1}=-\delta_{\mu_{0}%
}\varphi_{3}\\
\vdots\\
\varphi_{p}\circ\varphi_{p}+\sum_{\substack{1\leq i\leq p-1 \\}}\varphi
_{i}\circ\varphi_{2p-i}+\varphi_{2p-i}\circ\varphi_{i}=-\delta_{\mu_{0}%
}\varphi_{2p}\\
\sum_{1\leq i\leq p}\varphi_{i}\circ\varphi_{2p+1-i}+\varphi_{2p+1-i}%
\circ\varphi_{i}=-\delta_{\mu_{0}}\varphi_{2p+1}\\
\vdots
\end{array}
\right.  .
\]
Then the first term $\varphi_{1}$ of a deformation $\mu_{t}$ of a Lie algebra
law $\mu_{0}$ belongs to $Z^{2}(\mu_{0},\mu_{0}).$ This term is called the
infinitesimal part of the deformation $\mu_{t}$ of $\mu_{0}.$

\begin{definition}
A formal deformation of $\mu_{0}$ is called linear deformation if it is of
lenght one, that is of the type $\mu_{0}+t\varphi_{1}$ with $\varphi_{1}\in
Z^{2}(\mu_{0},\mu_{0}).$
\end{definition}

For a such deformation we have necessarily $\varphi_{1}\circ\varphi_{1}=0$
that is $\varphi_{1}\in L^{n}.$

Now consider $\varphi_{1}\in Z^{2}(\mu_{0},\mu_{0})$ for $\mu_{0}\in L^{n}.$
It is the infinitesimal part of a formal deformation of $\mu_{0}$ if and only
if there are $\varphi_{i}\in C^{2}(\mu_{0},\mu_{0}),$ $i\geq2 $, such that the
system $(I)$ is satisfied.

\begin{proposition}
If $H^{3}(\mu_{0},\mu_{0})=0$ then every $\varphi_{1}\in Z^{2}(\mu_{0},\mu
_{0})$ is an infinitesimal part of a formal deformation of $\mu_{0}.$
\end{proposition}

In fact, if $\varphi_{1}\in Z^{2}(\mu_{0},\mu_{0})$ then $\varphi_{1}%
\circ\varphi_{1}\in Z^{3}(\mu_{0},\mu_{0}).$ If $H^{3}(\mu_{0},\mu_{0})=0$,
then it exits $\varphi_{2}\in C^{2}(\mu_{0},\mu_{0})$ such that $\varphi
_{1}\circ\varphi_{1}=\delta\varphi_{2}.$ In this case $\varphi_{1}\circ
\varphi_{2}+\varphi_{2}\circ\varphi_{1}\in Z^{3}(\mu_{0},\mu_{0}).$ It exits
$\varphi_{3}\in C^{2}(\mu_{0},\mu_{0})$ such that
\[
\varphi_{1}\circ\varphi_{2}+\varphi_{2}\circ\varphi_{1}=\delta\varphi_{3}.
\]
As this we can solve step by step all the equations of the system $(I).$

Let us consider two formal deformations $\mu_{t}^{1}$ and $\mu_{t}^{2}$ of a
law $\mu_{0}.$ They are called equivalent if there exits a formal linear
isomorphism $\Phi_{t}$ of $\mathbb{C}^{n}$ of the following form
\[
\Phi_{t}=Id+\sum_{i\geq1}t^{i}g_{i}%
\]
with $g_{i}\in gl(n,\mathbb{C})$ such that
\[
\mu_{t}^{2}(X,Y)=\Phi_{t}^{-1}(\mu_{t}^{1}(\Phi_{t}(X),\Phi_{t}(Y))
\]
for all $X;Y\in\mathbb{C}^{n}.$

\begin{definition}
A deformation $\mu_{t}^{{}}$ of $\mu_{0}$ is called trivial if it is
equivalent to $\mu_{0}.$
\end{definition}

Let $\mu_{t}^{1}=\mu_{0}+\sum_{t=1}^{\infty}t^{i}\varphi_{i}$ and $\mu_{t}%
^{2}=\mu_{0}+\sum_{t=1}^{\infty}t^{i}\psi_{i}$ be two equivalent deformation
of $\mu_{0}.$ It is easy to see that
\[
\varphi_{1}-\psi_{1}\in B^{2}(\mu_{0},\mu_{0}).
\]
Thus we can consider that the set of infinitesimal parts of deformations is
parametrized by $H^{2}(\mu_{0},\mu_{0}).$

\subsection{Characteristic sequence of a nilpotent Lie algebra}

Let $\frak{n}$ be a complex finite dimensional Lie algebra. Consider the
derived subalgebra $C^{1}\frak{n}$. Let $Y\in\frak{n}-C^{1}\frak{n}$ be a vector of $\frak{n}$ which does not belong to the derived
subalgebra. Consider the ordered sequence
\[
c(Y)=(h_{1},h_{2},\cdots,)
\]
$h_{1}\geq h_{2},...,\geq h_{p},$ where $h_{i}$ is the dimension of the
$i^{th}$ Jordan bloc of the nilpotent operator $adY$. As $Y$ is necessary an
eigenvector of $adY$, then $h_{p}=1$. Let $Y_{1}$ and $Y_{2}$ be in
$\frak{n}-C^{1}\frak{n}$. Let be $c(Y_{1})=(h_{1},...,h_{p_{1}}=1)$
and $c(Y_{2})=(k_{1},...,k_{p_{2}}=1)$ the corresponding sequences. We have
$h_{1}\geq h_{2}\geq...\geq h_{p_{1}}$ and $k_{1}\geq k_{2}\geq...\geq
k_{p_{2}}$ with $h_{1}+...+h_{p_{1}}=k_{1}+...+k_{p_{2}}=n=\dim\frak{n}$. We
will say that $c(Y_{1})\geq c(Y_{2})$ if there is $i$ such that $h_{1}=k_{1},$
$h_{2}=k_{2},$ ... , $h_{i-1}=k_{i-1},$ $h_{i}>k_{i}.$ This defines a total
order relation on the set of sequences $c(Y)$ (lexicografic order) and we can
consider the maximum of these sequences.

\begin{definition}
The characteristic sequence of the nilpotent Lie algebra $\frak{n}$ is the
following sequence :
\[
c(\frak{n})=Sup\left\{c(Y),Y\in\frak{n}-C^{1}\frak{n}\right\}
\]
\end{definition}

It is an invariant up to isomorphism of $\frak{n}$, fintroduced by Ancochea and Goze in [5]. A vector $X\in
\frak{n}-C^{1}\frak{n}$ such that $c(X)=c(\frak{n})$ is called a
characteristic vector of $\frak{n}$.

This invariant is well adapted for study the deformations of nilpotent Lie
algebras.\ In fact let $\frak{n}$ and $\frak{n}{^{\prime}}$ be two
$n$-dimensional complex nilpotent Lie algebras and $\mu$ and $\mu^{\prime}$
the corresponding laws. Suppose that $\mu^{\prime}$ is a perturbation of $\mu
$, that is, in a fixed basis, the structural constant of $\mu^{\prime}$ are
close of those of $\mu$. In this case, the linear operator $ad_{\mu^{\prime}%
}Y$ is a perturbation (in the classical sense) of the linear operator
$ad_{\mu}Y$. As these two operators are nilpotent, the restriction of
$ad_{\mu^{\prime}}Y$ to the first Jordan block $J_{h_{1}}$of $ad_{\mu}Y $
satisfies $(ad_{\mu^{\prime}}Y\mid J_{h_{1}})^{h_{1}-2}\neq0$. Then, the first
Jordan block of $ad_{\mu^{\prime}}Y$ has a dimension greater or equal than
$h_{1}.$ This proves that
\[
c(\frak{n}{^{\prime}})\geq c(\frak{n}).
\]

\begin{proposition}
If $\frak{n}$ and $\frak{n}{^{\prime}}$ are two $n$-dimensional complex
nilpotent Lie algebras such that $\frak{n}^{\prime}$ is a perturbation of
$\frak{n}$, then
\[
c(\frak{n}{^{\prime}})\geq c(\frak{n}).
\]
\end{proposition}

This last property allowed to determine, for example, all the irreducible
components of the algebraic variety of $n$-dimensional nilpotent Lie algebras
for $n\leq8.$

\section{Characteristically nilpotent Lie algebras}

\bigskip 
In studying the varieties of laws, the characteristically nilpotent algebras
have shown their importance in the determination of irreducible components.
For example, in dimension $7$ there are two components, the first formed by
filiform Lie algebras and the second generated by the orbit closure of a
family of characteristically nilpotent Lie algebras [8].

The main problem in the study of characteristically nilpotent Lie algebras is
the determination of conditions for an algebra of derivations to be nilpotent:
for an arbitrary nilpotent Lie algebra the structure of the algebra of
derivations can variate from representations of the special linear algebras
$\frak{sl}_{n}$ to certain nilpotent Lie algebras.\newline The origin of all this is the cited result of Jacobson [48].

\begin{theorem}
Let $\frak{g}$ be a Lie algebra and suppose that it admits a nondegenerate derivation $f$. Then $\frak{g}$ is a nilpotent Lie algebra.
\end{theorem}

According to our convention, the Lie algebra is defined over a the field of complex numbers. Otherwise the assertion would be false, as it has been verified that this result fails when the characteristic of the base field is nonzero. \newline 
The example of Dixmier and Lister, appearing as the first known
characteristically nilpotent Lie algebra, was the response to the validity
question of Jacobson's theorem of 1955. This algebra is very interesting in
many aspects; it is one of the few known CNLA of nilindex 3, which is the
lowest possible nilindex such an algebra can have. We find this intriguing;
the authors not only gave the first example to a new class of nilpotent Lie
algebras, that also developed an ''extreme'' example in that sense.
Unfortunately, we do not know how Dixmier and Lister came to this
algebra.\newline The construction is of an eight dimensional Lie algebra
$\frak{g}_{0}$ defined over the basis $\left\{  X_{1},..,X_{8}\right\}  $ and law%
\begin{align*}
\left[  X_{1},X_{2}\right]   &  =X_{5};\;\left[  X_{1},X_{3}\right]
=X_{6};\;\left[  X_{1},X_{4}\right]  =X_{7};\;\left[  X_{1},X_{5}\right]
=-X_{8};\\
\left[  X_{2},X_{3}\right]   &  =X_{8};\;\left[  X_{2},X_{4}\right]
=X_{6};\;\left[  X_{2},X_{6}\right]  =-X_{7};\;\left[  X_{3},X_{4}\right]
=-X_{5};\\
\left[  X_{3},X_{5}\right]   &  =-X_{7};\;\left[  X_{4},X_{6}\right]  =-X_{8}.
\end{align*}

Let us define the following generalization of the central descending sequence
for a Lie algebra $\frak{g}$ :%
\[
\frak{g}^{[1]}=Der\left(  \frak{g}\right)  \left(  \frak{g}\right)  =\left\{
X\in\frak{g}\;|\;X=f\left(  Y\right)  ,\;f\in Der\left(  \frak{g}\right)
,\;Y\in\frak{g}\right\}
\]
and
\[
\frak{g}^{[k]}=Der\left(  \frak{g}\right)  \left(  \frak{g}^{[k-1]}\right)
,\;k>1
\]
The main result about this algebra is the following

\begin{theorem}
If $f$ is a derivation of $\frak{g}_{0}$ then $f\left(  \frak{g}_{0}\right)  \subset
C^{1}\frak{g}_{0}$; hence any derivation if nilpotent.
\end{theorem}

The proof of this is strongly related with the fact that the algebra $\frak{g}_{0}$ annihilates a power of the preceding sequence. For
this reason, they defined characteristically nilpotent Lie algebras as follows :

\begin{definition}
A Lie algebra $\frak{g}$ is called characteristically nilpotent if there
exists an integer $m$ such that $\frak{g}^{[m]}=0$.
\end{definition}

The listed algebra has a twelve dimensional Lie algebra of derivations, from
which six correspond to inner derivations. Now, any linear operator sending
the algebra $\frak{g}_{0}$ into its center, which is generated by the vectors
$X_{7}$ and $X_{8}$, is easily seen to be a derivation of $\frak{g}_{0}$. The
ideal uf these derivations has dimension eight, having a two dimensional
subspace in common with the space of inner derivations. This fact can be
interpreted in the sense that the dimension of the cohomology space
$H^{1}\left(  \frak{g}_{0},\frak{g}_{0}\right)  $ is as small as possible. Dixmier and
Lister asked if the algebras of this type, which satisfy the generalization of
the central descending sequence, were more treatable than ordinary nilpotent
Lie algebras. In certain aspects this is true, as the topological properties
of CNLA show, but on the other their determination and classification is a
rather difficult question, and one can hardly say that it constitutes a
simplification. However, CNLA have undoubtly contributed to a better
understanding of the geometry of the variety $\mathcal{N}^{n}$. The theorem
above proves in fact much more than the characteristic nilpotence of the
listed algebra : $\frak{g}_{0}$ is not the derived subalgebra of any Lie algebra.
Thus one can pose the question : if $\frak{g}$ is a CNLA, is it true that
$\frak{g}$ cannot be the derived subalgebra of a Lie algebra ? A first
condition is given above, as the nilpotence of $\frak{g}$ and the fact that
any derivation maps the algebra into its derived subalgebra ensures that it
cannot be a commutator algebra. Leger and T\^{o}g\^{o} found out other
conditions to assure the nonexistence of an algebra containing a given CNLA as
derived subalgebra [66] :

\begin{proposition}
Let $\frak{g}$ be a CNLA. If $Der\left(  \frak{g}\right)  $ annihilates the
center $Z\left(  \frak{g}\right)  $ of $\frak{g}$, then $\frak{g}$ is not a
derived algebra.
\end{proposition}

This proposition is based on the fact that for a CNLA we have
\[
\left[  \frak{g},Z_{i}\right]  \subset Z_{i-2}%
\]
where $Z_{0}=\left(  0\right)  $ and $Z_{i}$ is the largest subspace of
$\frak{g}$ such that $Der\left(  \frak{g}\right)  Z_{i}\subset Z_{i-1}$ for
$i\geq1$. The existence of an index such that $\frak{g}=Z_{r}$ follows
immediately. The authors also deduce an interesting numerical condition, also
based on this inductively defined sequence :

\begin{theorem}
Let $\frak{g}$ be a CNLA, and $n$ and $m$ be the smallest integers for which
$C^{n-1}\frak{g}=0$ and $\frak{g}^{[m]}=0$. If $2\left(  m-1\right)  >n+1$,
then $\frak{g}$ is not a derived subalgebra.
\end{theorem}

In particular, it follows that $\frak{g}$ is no derived subalgebra if
$Der\left(  \frak{g}\right)  \frak{g}\subset C^{1}\frak{g}$, which recovers
the property of Dixmier and Lister's algebra, or $\frak{g}^{[4]}=0$.
\newline Now, for the general case E. Luks [69] proved in 1976 that a CNLA
can appear as derived subalgebra of a Lie algebra. This remarkable fact
divides in fact the CNLA into two classes : those being commutators of others
and those not. The algebra $L$ of Luks has dimension $16$, and is defined over
the basis $\left\{  X_{1},..,X_{16}\right\}  $ by the law :%
\begin{align*}
\left[  X_{1},X_{i}\right]   &  =X_{5+i},\;i=2,3,4,5;\;\left[  X_{1}%
,X_{6}\right]  =X_{13};\;\left[  X_{1},X_{i}\right]  =X_{8+i},\;i=7,8\\
\left[  X_{2},X_{i}\right]   &  =X_{i+8},\;i=3,4,6;\;\left[  X_{2}%
,X_{5}\right]  =X_{15};\;\left[  X_{2},X_{7}\right]  =-X_{16}\\
\left[  X_{3},X_{4}\right]   &  =-X_{13}-\frac{9}{5}X_{15};\;\left[
X_{3},X_{5}\right]  =-X_{14};\;\left[  X_{3},X_{6}\right]  =-X_{16};\;\left[
X_{4},X_{5}\right]  =2X_{16}%
\end{align*}
The procedure used is the following : consider the transporter of an ideal $I$
into $J$, where both ideals are characteristic. Then the following statements
are verified

\begin{enumerate}
\item $I_{3}$ is the transporter of $C^{1}L$ to $0$.

\item $I_{4}$ is the transporter of $L$ to $Z\left(  L\right)  .$

\item $I_{2}$ is the transporter of $C^{1}L$ to $\left[  I_{4},I_{4}\right]  .$

\item $I_{6}$ is the transporter of $L$ to $\left[  I_{2},I_{2}\right]  .$

\item $I_{5}$ is the transporter of $I_{2}$ to $[L,I_{6}]$.
\end{enumerate}

where $I_{i}$ is the ideal generated by the vectors $\left\{  X_{j}\right\}
_{j\geq i}$ ( for $1\leq i\leq16$ ). Clearly $I_{2},..,I_{6}$ are
characteristic ideals, and it follows that $L$ is a CNLA. For the second part,
consider the derivations $f$ and $g$ defined respectively by :%
\begin{align*}
f\left(  X_{3}\right)   &  =X_{7};\;f\left(  X_{4}\right)  =2X_{8};\;f\left(
X_{5}\right)  =3X_{9}+2X_{11};\;f\left(  X_{6}\right)  =4X_{10}+5X_{12};\\
f\left(  X_{8}\right)   &  =X_{15};\;f\left(  X_{9}\right)  =2X_{16}%
;\;f\left(  X_{11}\right)  =-X_{16}.\\
g\left(  X_{i}\right)   &  =X_{i+1},\;i=2,3,4,5,7,11,13;\;g\left(
X_{8}\right)  =X_{9}+X_{11};\;g\left(  X_{9}\right)  =X_{10}+X_{12};\\
g\left(  X_{10}\right)   &  =X_{13}+X_{15};\;g\left(  X_{12}\right)
=-X_{13}-\frac{4}{5}X_{15};\;g\left(  X_{14}\right)  =-X_{16}%
\end{align*}
If one considers the brackets $\left[  f,g\right]  $ in $Der\left(  L\right)
$, it gives $ad\left(  X_{1}\right)  $ as result. Thus we can extend the
algebra to the semidirect product $\overline{L}=\left\{  f^{\prime},g^{\prime
}\right\}  +L$, where the brachets in $L$ are the same and the action of
$\left\{  f^{\prime},g^{\prime}\right\}  $ over $L$ is given by $\left[
f^{\prime},X_{i}\right]  =f\left(  X_{i}\right)  ,\;\left[  g^{\prime}%
,X_{i}\right]  =g\left(  X_{i}\right)  $ for all $i$ and $\left[  f^{\prime
},g^{\prime}\right]  =X_{1}$. It follows at once that $\left[  \overline
{L},\overline{L}\right]  =L$. \newline This algebra not only gives a
surprising response to the question of Dixmier and Lister, it moreover gives
an idea of how the algebra has to be structured for being candidate for
derived algebra. Observe that the key of the preceding cosntruction is the
existence of two derivations ( whose rest class modulo $B^{1}\left(
L,L\right)  $ is nonzero ) whose composition equals the adjoint operator of
the characteristic vector $X_{1}$. Thus, a necessary condition is deduced
immediately, namely the existence of derivations in the algebra such that
their composition is in the linear subspace generated by the adjoint operators
$ad\left(  X_{i}\right)  ,\;i=1,..,k$, where these vectors are generators of
the nilpotent Lie algebra. \newline Using a similar argumentation one can show
that Luks' algebra also has an unipotent automorphism group. One may ask if
there is a connection between the property of being a commutator algebra and
an automorphism group of this kind. L. Auslander remarked in [33] that
Dixmier and Lister's example has not an unipotent automorphism group. Now, a CNLA which is additionally a derived algebra and posseses a nonunipotent automorphism group could be constructed by considering the direct sum of two Lie algebras which satisfy the two first conditions [66].
Inspite of this result, there are wide
known classes of CNLA which cannot appear as a commutator algebra. This
concerns the filiform Lie algebras. It can be shown that if a filiform Lie
algebra $\frak{g}$ is the derived algebra of $\frak{g}^{\prime}$, then it
suffices to consider the case where $\dim\frak{g}^{\prime}=\dim\frak{g}+1$.
This has been done in [24]. The reduction is not difficult to prove, and
using it the assertion that $\frak{g}$ is not a derived algebra follows at
once. In fact, this reduction can be seen as a conseuence of a more general
result due to M. Goze and Y. B. Khakimdjanov [40], where they analyze in
detail the maximal tori of derivations of an arbitrary filiform Lie algebra. The following result can be interpreted as a characterization of those filiform Lie algebras which are not characteristically nilpotent : 

\begin{theorem}
Let $\frak{g}$ be an $\left(  n+1\right)  $-dimensional filiform Lie algebra
which has a nontrivial semisimple derivation $f$. There exists a basis
$\left\{  X_{0},..,X_{n}\right\}  $ adapted to $f$ such that the brackets of
$\frak{g}$ satisfy one of the following cases : 

\begin{enumerate}
\item $\frak{g}=L_{n}:$%
\[
\left[  X_{0},X_{i}\right]  =X_{i+1},\;1\leq i\leq n-1
\]

\item $\frak{g}=A_{n+1}^{r}\left(  \alpha_{1},..,\alpha_{t}\right)  \;1\leq
r\leq n-3,\;t=\left[  \frac{n-r-1}{2}\right]  $
\begin{align*}
\left[  X_{0},X_{i}\right]    & =X_{i+1},\;1\leq i\leq n-1\\
\left[  X_{i},X_{j}\right]    & =\left(  \sum_{k=1}^{t}\alpha_{k}\left(
-1\right)  ^{k-i}C_{j-k-1}^{k-1}\right)  X_{i+j+r},\;1\leq i,j\leq
n,\;i+j+r\leq n
\end{align*}

\item $\frak{g}=Q_{n}\;\ n=2m-1$%
\begin{align*}
\left[  X_{0},X_{i}\right]    & =X_{i+1},\;1\leq i\leq n-2\\
\left[  X_{i},X_{n-i}\right]    & =\left(  -1\right)  ^{i}X_{n},\;1\leq i\leq
n-1
\end{align*}

\item $\frak{g}=B_{n+1}^{r}\left(  \alpha_{1},..,\alpha_{t}\right)
,\;n=2m+1,\;1\leq r\leq n-4,\;t=\left[  \frac{n-r-2}{2}\right]  $%
\begin{align*}
\left[  X_{0},X_{i}\right]    & =X_{i+1},\;1\leq i\leq n-2\\
\left[  X_{i},X_{j}\right]    & =\left(  \sum_{k=1}^{t}\alpha_{k}\left(
-1\right)  ^{k-i}C_{j-k-1}^{k-1}\right)  X_{i+j+r},\;1\leq i,j\leq
n-1,\;i+j+r\leq n-1
\end{align*}

\item $\frak{g}=C_{n+1}\left(  \alpha_{1},..,\alpha_{t}\right)
,\;n=2m+1,\;t=m-1$%
\begin{align*}
\left[  X_{0},X_{i}\right]    & =X_{i+1},\;1\leq i\leq n-2\\
\left[  X_{i},X_{n-i}\right]    & =\left(  -1\right)  ^{i}X_{n},\;1\leq i\leq
n-1\\
\left[  X_{i},X_{n-i-2k}\right]    & =\left(  -1\right)  ^{i}\alpha_{k}%
X_{n},\;1\leq k\leq m-1,\;1\leq i\leq n-2k-1
\end{align*}
where $\left(  \alpha_{1},..,\alpha_{t}\right)  $ are parameters satisfying
the polynomial relations given by the Jacobi relations over this basis.
\end{enumerate}
\end{theorem}

\section{Structural properties of CNLA}

After the example of Dixmier and Lister in 1957, Leger and T\^{o}g\^{o}
iniciated the structural study of CNLA. Their paper [66] does not provide
additional examples, but it is of considerable significance for later work. At
first, they observe that the property of being characteristically nilpotent
does not depend on the ground field. More precisely : if the Lie algebra
$\frak{g}$ is characteristically nilpotent as $F$-algebra ( here it is not
necessary to suppose that it has characteristic zero ) and $K\backslash F$ is
a field extension, then $\frak{g}$ is also a CNLA as $K$-algebra. However, the
structural properties deduced by the authors are more important, as they give
an idea of which algebras have to be avoided in the search after CNLA :

\begin{lemma}
If $\frak{g}$ is characteristically nilpotent, then

\begin{enumerate}
\item  the center $Z\left(  \frak{g}\right)  $ of $\frak{g}$ is contained in
the derived subalgebra $C^{1}\frak{g}$.

\item $C^{2}\frak{g}\neq0$.
\end{enumerate}
\end{lemma}

The first condition makes reference to the nonexistence of direct summands in
$\frak{g}$ which constitute of central vectors. Thus the study of
characteristically nilpotent Lie algebras reduces to nonsplit nilpotent Lie
algebras. The second condition has a more important consequence : it tells
that for a Lie algebra being characteristically nilpotent, the nilindex must
be at least three ( observe that this is the index for the algebra of Dixmier
and Lister ). This fact is remarkable, as it shows the incompatibility of
being as nilpotent as possible ( as it occurs for the 2-step nilpotent or
metabelian Lie algebras ) and having all its derivations nilpotent. Metabelian
Lie algebras and their derivations have been deeply studied by Leger and Luks
[64], where they proved that its rank is always greater than one, the
equality given only under certain conditions. Recently Galitski and Timashev
[37] have used invariant theory to classify these algebras up to dimension
nine. The preceding lemma leads to the question wheter a CNLA can be a direct
sum. The following result is also from [66] :

\begin{lemma}
Let $\frak{g}$ be a nilpotent Lie algebra. If $\frak{g}$ is the direct sum of
two nontrivial ideals, one of which is central, then it posseses at least
nontrivial semisimple derivation.
\end{lemma}

These two lemmas give the following reinterpretation of the sequence
$\frak{g}^{[k]}$ introduced earlier :

\begin{theorem}
Let $\frak{g}$ be a Lie algebra and $Der\left(  \frak{g}\right)  $ its Lie
algebra of derivations. Then $\frak{g}$ is characteristically nilpotent if and
only if $Der\left(  \frak{g}\right)  $ is nilpotent and $\dim\frak{g}\geq2$.
\end{theorem}

It follows from the proof that if all derivations of $\frak{g}$ are nilpotent,
then $\frak{g}$ is also a nilpotent Lie algebra. Thus the characteristic
nilpotence is a phenomena which can only be observed in the variety of
nilpotent Lie algebra laws $\mathcal{N}^{n}$. The theorem can be reformulated
by saying that the holomorph $H\left(  \frak{g}\right)  $ of $\frak{g}$ is
nilpotent, where the holomorph is the split extension of $Der\left(
\frak{g}\right)  $ by $\frak{g}$. In connection with metabelian Lie algebras,
this reformulation says that for a 2-step nilpotent Lie algebra the holomorph
cannot be nilpotent. The holomorph is also useful to describe properties valid
also for solvable Lie algebras, as the following

\begin{theorem}
Let $\frak{g}$ be a Lie algebra. If a Cartan subalgebra $H$ of $\frak{g}$ is
characteristically nilpotent, then $\frak{g}$ is a solvable Lie algebra.
\end{theorem}

As noted by the authors, the algebra $\frak{g}$ can be solvable non-nilpotent.
We remark that this theorem has been generalized in 1961 by S. T\^{o}g\^{o}
[94].

It has often been asked wheter CNLA exist for any possible dimension. The
answer is in the affirmative, and in fact it was enough to find examples of
dimension $7\leq n\leq13$ to derive its existence in any dimension. The key
result was the possibility of a decomposition into smaller blocks that have
also the property of being characteristically nilpotent, as done in the
classical theory :

\begin{theorem}
Let $\frak{g}$ $=\bigoplus_{i=1}^{n}\frak{g}_{i}$ be a direct sum of ideals.
Then $\frak{g}$ is characteristically nilpotent if and only if $\frak{g}_{i}$
is characteristically nilpotent for $1\leq i\leq n$.
\end{theorem}

As said, having examples from dimensions seven to thirtheen, the direct sums
of them give CNLA in any dimension. The nine dimensional example was given by
J. Dyer in 1970, in connection with her study of nilpotent Lie groups which
have expanding automorphisms. Over the basis $\left\{  X_{1},..,X_{9}\right\}
$ the Lie algebra is given by
\begin{align*}
\left[  X_{1},X_{2}\right]    & =X_{3};\;\left[  X_{1},X_{3}\right]
=X_{4};\;\left[  X_{1},X_{5}\right]  =X_{7};\;\left[  X_{1},X_{8}\right]
=X_{9};\\
\left[  X_{2},X_{3}\right]    & =X_{5};\;\left[  X_{2},X_{4}\right]
=X_{7};\;\left[  X_{2},X_{5}\right]  =X_{6};\;\left[  X_{2},X_{7}\right]
=-X_{8};\\
\left[  X_{3},X_{7}\right]    & =-\left[  X_{4},X_{5}\right]  =X_{9}%
\end{align*}
This was the first given CNLA with an unipotent automorphism group. Two years
later G. Favre constructed a seven dimensional example with the same property. 
This example is one of the three filiform CNLA in dimension $7$: 
\begin{align*}
\left[  X_{1},X_{i}\right]    & =X_{i+1},\;2\leq i\leq6\\
\left[  X_{3},X_{2}\right]    & =X_{6}\\
\left[  X_{4},X_{2}\right]    & =\left[  X_{5},X_{2}\right]  =X_{7}\\
\left[  X_{4},X_{3}\right]    & =-X_{7}%
\end{align*}
To complete the construction of CNLA, there remains to find examples in dimensions $10-13$. These were given by Luks using computational methods
[68]. Once the question of their existence in any possible dimension, we can
ask even more : for any possible nilindex $p\geq3$, does there exist a CNLA in
any dimension? In [7] the question is answered in the affirmative for $p=5$.
This is a consequence of the classification of nilpotent Lie algebras of
characteristic sequence $\left(  5,1,..,1\right)  $ whose derived subalgebra
is non-abelian. In fact, we prove that if a Lie algebra $\frak{g}$  with this
characteristic sequence is characteristically nilpotent, then it satisfies
$D^{2}\frak{g}\neq0$. \newline In 1961 T\^{o}g\^{o} published a paper
reviewing most of known results about the derivation algebras of Lie algebras
( over a field of charateristic zero ). He also gives an example about two
nonisomorphic Lie algebras whose Lie algebra of derivations is the same,
proving in that manner that a Lie algebra is not entirely determined by its
derivations. Among various results about classical and reductive algebras, he
also generalizes the concept of CNLA to characteristically solvable Lie
algebras [94]. However, here we are only concerned with results about
nilpotent Lie algebras. An often asked question is the relation between a Lie
algebra $\frak{g}$ which is a ( finite ) sum of ideals and the structure of
$Der\left(  \frak{g}\right)  $. To this respect, in [94] the following
theorem is proved :

\begin{theorem}
Let $\frak{g}=\bigoplus_{i=1}^{n}\frak{g}_{i}$ be a direct sum of ideals. Then
$Der\left(  \frak{g}\right)  =\bigoplus_{i=1}^{n}Der\left(  \frak{g}%
_{i}\right)  $ if and only if $\frak{g}$ satisfies one of the following
conditions :

\begin{enumerate}
\item $Z\left(  \frak{g}\right)  =\left(  0\right)  $

\item $\frak{g}$ is a perfect Lie algebra ( i.e. $\frak{g}=\left[
\frak{g},\frak{g}\right]  $ )

\item  All the $\frak{g}_{i}$'s except one is such that $Z\left(  \frak{g}%
_{i}\right)  =\left(  0\right)  $ and $\frak{g}_{i}=\left[  \frak{g}%
_{i},\frak{g}_{i}\right]  $.
\end{enumerate}
\end{theorem}

For a nilpotent Lie algebra $\frak{g}$, this implies that the structure of its
derivations is more than the sum of the derivations corresponding to its
summands. The following proposition gives the precise form of $Der\left(
\frak{g}\right)  $ :

\begin{proposition}
Let $\frak{g}=\bigoplus_{i=1}^{n}\frak{g}_{i}$ be a direct sum of ideals.
Then
\[
Der\left(  \frak{g}\right)  =\bigoplus_{i=1}^{n}\left(  Der\left(
\frak{g}_{i}\right)  \oplus\left(  \bigoplus_{i\neq j}\mathcal{D}\left(
\frak{g}_{i},\frak{g}_{j}\right)  \right)  \right)
\]
where
\[
\mathcal{D}\left(  \frak{g}_{i},\frak{g}_{j}\right)  =\left\{  h\in End\left(
\frak{g}\right)  \;|\;h\left(  \frak{g}_{k}\right)  =0\text{ if }k\neq
i,h\left(  \frak{g}_{i}\right)  \subset Z\left(  \frak{g}_{i}\right)  \text{
and }h\left(  \left[  \frak{g}_{i},\frak{g}_{j}\right]  \right)  =0\right\}
\]
\end{proposition}

Thus if one of the conditions in theorem 8 is satisfied, then $\mathcal{D}%
\left(  \frak{g}_{i},\frak{g}_{j}\right)  $ vanishes. \newline In the same
paper T\^{o}g\^{o} presents a list of problems of interest, specially in
connection with CNLA : do there exists CNLA of derivations? From the structure of derivations for the example of Dixmier and Lister, as well as the scarceness of outer derivations, it is obvious that this algebra does not have a characteristically nilpotent algebra of derivations. As to our knowledge, nobody has answered explicitely to this question until now, though
the answer is in the affirmative. In [8] we construct examples of CNLA of
derivations and generalize the question to higher indexes.

\begin{example}
Let $\frak{g}$ be the Lie algebra with associated law 
\begin{align*}
\mu_{5}\left(  X_{1},X_{i}\right)  =X_{i+1},\;i\in\{2,3,4,5\};\;\mu
_{5}\left(  X_{5},X_{2}\right)  =\mu_{5}\left(  X_{3},X_{4}\right)  =X_{6} \\
\mu_{5}\left(  X_{7},X_{3}\right)  =X_{6}\;\ \mu_{5}\left(X_{7},X_{2}\right)  =X_{5}+X_{6}.
\end{align*}
over the basis $\left\{X_{1},..,X_{7}\right\}$. The Lie algebra of derivations  
$Der\left(  \frak{g}\right)  $ is ten dimensional and isomorphic to 
\[
\text{%
\begin{tabular}
[c]{lll}%
$\lbrack Z_{1},Z_{2}]=Z_{3},$ & $[Z_{2},Z_{6}]=-Z_{5},$ & $\left[  Z_{7}%
,Z_{8}\right]  =2Z_{5}-2Z_{6}+2Z_{10}$\\
$\left[  Z_{1},Z_{3}\right]  =Z_{4},$ & $[Z_{2},Z_{8}]=-Z_{6},$ &
$[Z_{7},Z_{9}]=Z_{5}-2Z_{6}+2Z_{10}$\\
$\lbrack Z_{1},Z_{4}]=Z_{5},$ & $[Z_{2},Z_{9}]=-Z_{4}-2Z_{6},$ & $[Z_{8}%
,Z_{9}]=2Z_{6}-2Z_{10}$\\
$\lbrack Z_{1},Z_{7}]=-Z_{4},$ & $[Z_{2},Z_{10}]=-Z_{5},$ & \\
$\lbrack Z_{1},Z_{8}]=-Z_{6},$ & $[Z_{3},Z_{8}]=-Z_{5},$ & \\
& $[Z_{3},Z_{9}]=-Z_{5},$ &
\end{tabular}
}%
\]
It is routine to verify that this algebra is a CNLA.
\end{example}

Among many other examples, we present the following, which is important in connection 
with the study of irreducible components of the variety $\frak{N}^{n}$ :

\begin{theorem}
For any  $\alpha\in\mathbb{C}-\{0,2\}$  the family of nilpotent Lie algebras
given by 
\begin{align*}
\lbrack X_{1},X_{i}]  & =X_{i+1},\;2\leq i\leq5\\
\lbrack X_{4},X_{2}]  & =\alpha X_{6};\;\\
\lbrack X_{3},X_{2}]  & =\alpha X_{5}+X_{7}\\
\lbrack X_{7},X_{3}]  & =X_{6}\\
\lbrack X_{7},X_{2}]  & =X_{5}+X_{6}%
\end{align*}
has a characteristically nilpotent Lie algebra of derivations. 
\end{theorem}
This follows at once from the fact that the derivatiosn are given by :
\[%
\begin{array}
[c]{lll}%
\left[  Z_{1},Z_{2}\right]  =Z_{3}, & \left[  Z_{2},Z_{3}\right]  =-\alpha
Z_{5}-Z_{6}, & \left[  Z_{7},Z_{9}\right]  =2Z_{5},\\
\left[  Z_{1},Z_{3}\right]  =Z_{4}, & \left[  Z_{2},Z_{6}\right]  =-Z_{5}, &
\left[  Z_{7},Z_{10}\right]  =Z_{5},\\
\left[  Z_{1},Z_{4}\right]  =Z_{5}, & \left[  Z_{2},Z_{10}\right]
=-Z_{4}-\alpha Z_{5}, & \left[  Z_{9},Z_{10}\right]  =\frac{2}{\alpha}%
Z_{8}+\frac{2}{\alpha}Z_{4},\\
\left[  Z_{1},Z_{7}\right]  =-Z_{4}, & \left[  Z_{3},Z_{9}\right]  =-Z_{5}, &
\\
\left[  Z_{1},Z_{8}\right]  =-Z_{5}, & \left[  Z_{3},Z_{10}\right]  =-Z_{5}, &
\\
\left[  Z_{1},Z_{9}\right]  =-Z_{6}, &  & \\
\left[  Z_{1},Z_{10}\right]  =-Z_{8}-Z_{4}, &  &
\end{array}
\]

This examples, as well as other considered in [8] have a common property : there always exists an outer derivation $\theta$ which belongs to the derived subalgebra of $Der\left(\frak{g}\right)$. This and the method used to deduce the examples have led to the 
 
\begin{conjecture}
If $\frak{g}$ is a CNLA of derivations, then there exist outer derivations 
$\theta_{1},\theta_{2},\theta_{3}$ such that
\[
\lbrack\theta_{1},\theta_{2}]=\lambda\theta_{3}\;\;\left(  \operatorname{mod}%
\;IDer\left(  \frak{g}\right)  \right)
\]
where $\lambda\in\mathbb{C}-\{0\}$ e $\ IDer\left(  \frak{g}\right)  $ denotes
the space of inner derivations.
\end{conjecture}

We now come to the generalization announced. Let $Der^{[k]}\frak{g}=Der\left(  Der\left(  ...Der\left(  \frak{g}\right)
\right)  \right)  $ be the  $k$-th Lie algebra of derivations. Thus we have the sequence
\[
\left(  Der\left(  \frak{g}\right)  ,Der^{[2]}\frak{g},....,Der^{[k]}%
\frak{g},....\right)
\]

\begin{definition}
A Lie algebra $\frak{g}$ is called characteristically nilpotent of index $k$
if the $\left(  k-1\right)  ^{th}$ Lie algebra of derivations $Der^{[k-1]}%
\frak{g}$ is characteristically nilpotent.
\end{definition}

\begin{remark}
It would be of great interest to know if there exist CNLA of infinite index, as this would us give the possibility to develop a theory analogue to Schenkman's one [85] for these algebras. 
The structure of the variety of filiform Lie algebras $\mathcal{F}^{m}$ for $m\geq8$
seems to suggest the existence of such algebras, but there is no manner to prove it. Observe that the determination of such an algebra is far from being a computational problem. The question is more to find a new invariant which measures which is the gratest possible index, if any. Up to the moment, the biggest index known is 5.
\end{remark}

\section{Subspaces of CNLA}

Around 1984, when some authors had already constructed infinite families of
CNLA, the interest on these algebras turned to its topological and geometrical
properties. R. Carles proved in [19] the following result :

\begin{proposition}
The CNLA constitute a constructible set of the variety $\frak{N}^{n}$ which is
empty for $n\leq6$ and nonempty for $n\geq7$. 
\end{proposition}

This proposition is another way to prove the existence of CNLA in arbitrary
dimension, and its advantage is being independent from any example. Its proof
is based on the conjugacy classes of maximal tori of derivations over a
nilpotent Lie algebra of dimension $n$ ( [36] ), as well as the action of the
general linear group $GL\left(  n,\mathbb{C}\right)  $ on $\frak{g}$  ( the
result is in fact true for any algebraically closed field of characteristic
zero ). The seven dimensional CNLA given by Favre in 1972 is generalized
in the following manner : over the basis $\left\{  X_{1},..,X_{n}%
,X_{n+1}\right\}  $ the Lie algebra structure is given by :%
\begin{align*}
\left[  X_{1},X_{i}\right]    & =X_{i+1},\;2\leq i\leq n\\
\left[  X_{2},X_{3}\right]    & =X_{n}\\
\left[  X_{2},X_{4}\right]    & =\left[  X_{2},X_{5}\right]  =-\left[
X_{3},X_{4}\right]  =X_{n+1}%
\end{align*}
For $n=6$ Favre's example is recovered. The interest of this family is that it
is obtained by considering central extensions of an algebra $\frak{g}^{\prime
}$ by $\mathbb{C}$, which proves the power of extension theory for the study
of CNLA. It is also proven that any extension by the center of a CNLA is also
characteristically nilpotent, where an extension by the center is a central
extension of a Lie algebra $\frak{g}$ by $\mathbb{C}^{p}$ whose center is
isomorphic to $\mathbb{C}^{p}$. The same procedure has been used in [22] to
obtain lots of CNLA in arbitrary dimension and mixed characteristic
sequences. Carles also remarks that the set of CNLA is never closed, which is
immediate from the preceding, and for the particular case of dimension $7$ he
proves that it is neither open. In [8] we have extended this result to any
dimension : 

\begin{theorem}
For $n\geq8$ the set $\mathcal{S}_{n}$ of CNLA is not open in the variety
$\frak{N}^{n}$.
\end{theorem}

The family constructed is based on the results of the classification of 8
dimensional filiform Lie algebra due to Goze and the first author [4]. Also
the deformation structure is based on this result :

Let $\frak{g}_{n,17}$ $\left(  n\geq8\right)  $ be the Lie algebra defined by
the brackets
\begin{align*}
\lbrack X_{1},X_{i}] &  =X_{i+1},\;1\leq i\leq n-1\\
\lbrack X_{4},X_{2}] &  =X_{n},\\
\lbrack X_{3},X_{2}] &  =X_{n-1}+X_{n}%
\end{align*}
It is immediate that the algebra is filiform and characteristically nilpotent.
Let $\psi\in Z^{2}\left(  \frak{g}_{n,17},\frak{g}_{n,17}\right)  $ be the
linear expandable cocycle defined by
\begin{align*}
\psi\left(  X_{5},X_{3}\right)   &  =X_{n},\;\psi\left(  X_{5},X_{2}\right)
=\psi\left(  X_{4},X_{3}\right)  =X_{n-1},\\
\psi\left(  X_{k},X_{2}\right)   &  =2X_{n-4+\left[  \frac{k}{2}\right]
},\;k=3,4
\end{align*}
\newline Let $\frak{g}_{n,17}+\varepsilon\psi$ be an infinitesimal deformation
of $\frak{g}_{n,17}.\;\;\;$\newline Now we consider the change of Jordan basis
$X_{1}^{\prime}=X_{1},\;X_{2}^{\prime}=X_{2}+a_{3}X_{3}+a_{4}X_{4}+a_{5}X_{5}$
with the relations
\begin{align*}
1+a_{3}^{2}\varepsilon-2\varepsilon a_{4} &  =0\\
3a_{5}\varepsilon+a_{3}a_{4}-a_{3}^{2}\varepsilon-2\varepsilon a_{4} &  =0
\end{align*}
Written in the new basis the algebra $\frak{g}_{n,17}+\varepsilon\psi$ is
isomorphic to the Lie algebra $\frak{g}_{n,18}$ defined by
\begin{align*}
\lbrack X_{1},X_{i}] &  =X_{i+1},\;1\leq i\leq n-1\\
\lbrack X_{5},X_{3}] &  =\varepsilon X_{n}\\
\lbrack X_{5},X_{2}] &  =[X_{4},X_{3}]=\varepsilon X_{n-1}\\
\lbrack X_{4},X_{2}] &  =2\varepsilon X_{n-2}\\
\lbrack X_{3},X_{2}] &  =2\varepsilon X_{n-3}%
\end{align*}
From the linear system $\left(  S\right)  $ associated to this algebra [6]
we deduce the existence of nonzero eigenvalues for diagonalizable derivations
of $\frak{g}_{n,18},$ so it cannot be characteristically nilpotent.
\newline Following with the seven dimensional case, in [18] the irreducible
components of the variety $\frak{N}^{7}$ are analyzed in relation with
characteristic nilpotence. It is well known that this variety has two
irreducible components, one corresponding to the filiform Lie algebras,
$\mathcal{F}_{7}$, and one consisting of non-filiform Lie algebras. The
filiform components has only three CNLA, which obviously don't constitute a
dense subset as none of them has an open orbit, while for the other component there exists a family of CNLA
constituting a nonempty Zariski open subset. The family is precisely the one
given as example above. For $n\geq8$, the situation for $\mathcal{F}_{8}$
changes radically :

\begin{theorem}
Let $n\geq8$. Then any irreducible component $C$ of $\mathcal{F}_{n}$ contains
a nonempty Zariski open subset $A$ cosisting of CNLAs.
\end{theorem}

The result is proven in [40], where even more is obtained, namely that for
any open set in $\mathcal{F}_{n}$ a CNLA belonging to this set can be found.
Other versions relative to this have been treated by H. Kraft and Ch.
Riedtmann in [61]. Is it true that for dimensions $n\geq8$ any irreducible
component of the variety $\frak{N}^{n}$ contains an open subset of CNLAs. For
$n=8$ the response is affirmative, and can be found in [8].

\begin{proposition}
For any irreducible component $C$ of the variety $\frak{N}^{8}$ there exists a
nonempty Zariski-open subsets consisting of CNLAs.
\end{proposition}

As commented above, the 1-parameter family that defines the second irreducible
component of $\frak{N}^{7}$ has the property of being characteristically
nilpotent of index 2, according to the definition given. This leads naturally
to the question wheter there exist irreducible components in $\frak{N}^{n}$
$\left(  n\geq9\right)  $ which admit nonempty open susbsets formed by CNLA od
derivations. We finally remark that this problem is related with the tower
problem in group theory.  

\subsection{Characteristically nilpotent Lie algebras obtained from
nilradicals of Borel subalgebras}

\bigskip As commented earlier, the difficulty of constructing and
characterizing CNLA led many authors to conclude that they were scarce within
the variety of nilpotent Lie algebra laws, though results like those of Carles
[18] pointed out their importance. The question was definitively solved by
Khakimdjanov in 1988, in a series of papers ( [53], [54] ), where he
treated with the cohomology of parabolic subalgebras of simple Lie algebras,
first studied by Kostant in 1963 ( [58], [59] ), and applied these results
to the study of deformations of the nilradicals of Borel subalgebras of simple
Lie algebras. For classical topics we refer the reader to [46], [44], [23]
and [58]. \newline In [53] the author developes the cohomological tools
needed, such as the fundamental cohomologies, as well as adequate filtrations
for these spaces. In [54] this information is applied to prove that almost
all deformations of the cited nilradicals are CNLA. \newline Following the
notation used in [44], let $L$ be a simple Lie algebra of rank $l>1,\frak{\ }%
H$ its Cartan subalgebra, $\Phi$ the root system associated to $H$, $\Phi^{+}$
the system of positive roots relative to a certain ordering and $\Delta$ the
system of simple roots. Recall that a Borel subalgebra is a maximal solvable
subalgebra of $L$.

We consider the subalgebra $B\left(  \Delta\right)  =H+\coprod_{\alpha\in
\Phi^{+}}L_{\alpha}$ , where $L_{\alpha}$ is the root space corresponding to
the root $\alpha.$ This subalgebra is a Borel subalgebra of $L$ called
standard relative to the Cartan subalgebra $H$. Now any Borel subalgebra of
$L$ is conjugated to a standard Borel subalgebra [11], and if $\frak{n}$
denotes the nilradical of an algebra $\frak{g}$ we have $\frak{n}\left(
B\left(  \Delta\right)  \right)  =\coprod_{\alpha\in\Phi^{+}}L_{\alpha}.$
Define $\Phi\left(  i\right)  $ as
\[
\Phi\left(  i\right)  =\left\{  \alpha\in\Phi^{+}\;|\;\alpha=\alpha_{j_{1}%
}+..+\alpha_{j_{i}},\;\alpha_{j_{t}}\in\Delta\text{ for }1\leq t\leq
i\right\}
\]
Then we can define a graduation on $\frak{n}\left(  B\left(  \Delta D\right)
\right)  $\ by setting $F_{k}\frak{n}\left(  B\left(  \Delta\right)  \right)
=\bigoplus_{i\geq k}\frak{n}_{i}\left(  B\left(  \Delta\right)  \right)  ,$
where $\frak{n}_{i}\left(  B\left(  \Delta\right)  \right)  =\coprod
_{\alpha\in\Phi\left(  i\right)  }L_{\alpha}.$ The filtration in the space of
cochains is given by%
\[
F_{k}C^{j}\left(  \frak{n},\frak{n}\right)  =\left\{  c\in C^{j}%
(\frak{n},\frak{n)}\;|\;c\left(  a_{1},..,a_{j}\right)  \in F_{t_{1}%
+..+t_{j}+k}\frak{n}\right\}
\]
whenever $a_{i}\in F_{t_{i}}\frak{n}\left(  B\left(  \Delta\right)  \right)  $
and where $\frak{n}=\frak{n}\left(  B\left(  \Delta\right)  \right)  .$

This filtration extends to the cocycles and coboundaries.\ Until now we
exclude $L$ \ to be a simple algebra of the following types%
\[
A_{i}\;\left(  1\leq i\leq5\right)  ,\;B_{2},B_{3},C_{3},C_{4},D_{4},G_{2}%
\]
\newline The reason is that for these algebras certain identities among the
fundamental cohomologies and the spaces $F_{k}H^{j}\left(  \frak{n}%
,\frak{n}\right)  $ for $k=0,1$ do not coincide [53, theorem 6]. For example,
for those algebras excluded and distinct from $A_{i}\;\left(  i=1,2,3\right)
,G_{2}$ the cohomology space $F_{0}H^{2}\left(  \frak{n},\frak{n}\right)  $ is
not zero. On the other side, it is shown that the following system of cocycles
suffices for a set of representatives of a basis of $F_{0}H^{2}\left(
\frak{n},\frak{n}\right)  :\;\left\{  f_{\alpha,\beta}\;|\;\left(
\alpha,\beta\right)  \in E\right\}  $ with
\[
f_{\alpha,\beta}\left(  x_{\gamma},x_{\delta}\right)  =\left\{
\begin{tabular}
[c]{ll}%
$x_{\sigma_{\alpha}\sigma_{\beta}\left(  \delta\right)  }$ & for $\left(
\gamma,\delta\right)  =\left(  \alpha,\sigma_{\alpha}\left(  \beta\right)
\right)  $\\
$0$ & otherwise
\end{tabular}
\right.
\]
where $\sigma_{\alpha}$ is the involution associated to the root $\alpha$ and
$E$ is the set of pairs of simple roots $\left(  \alpha,\beta\right)  $ in
which $\left(  \alpha,\beta\right)  $ is identified with $\left(  \beta
,\alpha\right)  $ if $\alpha$ is not joined to $\beta$ in the Dynkin diagram.

\begin{theorem}
Let $L$ be a simple Lie algebra and $\frak{n}$ be the nilradical of a Borel
subalgebra. Let $\psi=\sum_{\omega\in E}\lambda_{\omega}f_{\omega}$ an element
of $F_{0}H^{2}\left(  \frak{n},\frak{n}\right)  $ with $\lambda_{\omega}\neq0$
for all $\omega.$ Then the Lie algebra $\frak{n}\left(  \psi\right)  $
obtained from the linearly expandable cocycle $\psi$ is characteristically nilpotent.
\end{theorem}

Let $L\in\left\{  A_{4},A_{5},B_{3},C_{3},C_{4},D_{4}\right\}  $. For these
algebras we have $F_{0}H^{2}\left(  \frak{n},\frak{n}\right)  \neq0$. In [53]
it is proven that the basis is composed by cocycles of the form $x_{\alpha}\wedge
x_{\beta}\rightarrow x_{\gamma}$, where $\alpha,\beta$ and $\gamma$ are
enumerated in the following table :%
\[%
\begin{tabular}
[c]{|l|l|l|l|}\hline
$L$ & $\alpha$ & $\beta$ & $\gamma$\\\hline
$A_{4}$ & $%
\begin{array}
[c]{c}%
\alpha_{1}\\
\alpha_{1}\\
\alpha_{4}%
\end{array}
$ & $%
\begin{array}
[c]{c}%
\alpha_{4}\\
\alpha_{1}+\alpha_{2}\\
\alpha_{3}+\alpha_{4}%
\end{array}
$ & $%
\begin{array}
[c]{c}%
\alpha_{2}+\alpha_{3}\\
\alpha_{2}+\alpha_{3}+\alpha_{4}\\
\alpha_{1}+\alpha_{2}+\alpha_{3}%
\end{array}
$\\\hline
$A_{5}$ & $%
\begin{array}
[c]{c}%
\alpha_{2}\\
\alpha_{4}%
\end{array}
$ & $%
\begin{array}
[c]{c}%
\alpha_{1}+\alpha_{2}\\
\alpha_{4}+\alpha_{4}%
\end{array}
$ & $%
\begin{array}
[c]{c}%
\alpha_{3}+\alpha_{4}+\alpha_{5}\\
\alpha_{1}+\alpha_{2}+\alpha_{3}%
\end{array}
$\\\hline
$B_{3}$ & $\alpha_{1}$ & $\alpha_{1}+\alpha_{2}$ & $\alpha_{2}+2\alpha_{3}%
$\\\hline
$C_{3}$ & $\alpha_{1}$ & $\alpha_{1}+\alpha_{2}$ & $2\alpha_{2}+\alpha_{3}%
$\\\hline
$C_{4}$ & $\alpha_{2}$ & $\alpha_{1}+\alpha_{2}$ & $2\alpha_{3}+\alpha_{4}%
$\\\hline
$D_{4}$ & $%
\begin{array}
[c]{c}%
\alpha_{1}\\
\alpha_{3}\\
\alpha_{4}%
\end{array}
$ & $%
\begin{array}
[c]{c}%
\alpha_{1}+\alpha_{2}\\
\alpha_{2}+\alpha_{3}\\
\alpha_{2}+\alpha_{4}%
\end{array}
$ & $%
\begin{array}
[c]{c}%
\alpha_{2}+2\alpha_{3}\\
\alpha_{1}+\alpha_{2}+\alpha_{4}\\
\alpha_{1}+\alpha_{2}+\alpha_{3}%
\end{array}
$\\\hline
\end{tabular}
\]

\begin{theorem}
Let $L$ be a simple Lie algebra of types $A_{4},A_{5},B_{3},C_{3},C_{4}$ or
$D_{4}$. Let $\frak{n}$ be the nilradical of the standard Borel subalgebra
$B\left(  \Delta\right)  $ and $\varphi=\sum_{\omega\in E}\lambda_{\omega
}f_{\omega},$ where $\left\{  f_{\omega}\;|\;\omega\in E\right\}  $ is a basis
of $F_{0}H^{2}\left(  \frak{n},\frak{n}\right)  $ from the previous table,
with $\lambda_{\omega}\neq0$ for all $\omega\in E$. Then the nilpotent Lie
algebra defined by a deformation%
\[
\left[  X,Y\right]  _{t}=\left[  X,Y\right]  +t\varphi\left(  X,Y\right)
+t^{2}\varphi_{2}\left(  X,Y\right)  ,\;t\neq0
\]
is a CNLA.
\end{theorem}

These results are certainly of interest for the theory of CNLA. It provides not only a relation between the classical Cartan theory of Lie algebras, it moreover gives, in a certain manner, a natural interpretation of the characteristic nilpotence. On the other side, the frequency of CNLA in $\frak{N}^{n}$ is proven in an elegant manner. 

\section{Characteristically nilpotent filiform Lie algebras}

Most constructions of CNLA made are based on the deformation theory of the
naturally graded filiform Lie algebra $L_{n}$. The reason is not only its
simplicity; it turns out to have the most elementary law among the filiform
Lie algebras. Vergne proved in [96] that any filiform Lie algebra can be
obtained by a deformation of this algebra. For this reason this algebra has
been the preferred starting point for constructing families of CNLA [100],
[54], though recently other authors have turned their interest into the
deformations of the other naturally graded filiform Lie algebra [22].
\newline Certain results about the cohomologies of filiform Lie algebras are
contained in Vergne's paper [97]. Recall the notations introduced for the
filtered cohomology :

\begin{lemma}
Let $\frak{g}$ be a p-step nilpotent Lie algebra and  $d_{i}=\dim
\,F_{i}\frak{g}$. 

\begin{enumerate}
\item  If $j>d_{1}$, then $F_{r}Z^{j}\left(  \frak{g},\frak{g}\right)
=Z^{j}\left(  \frak{g},\frak{g}\right)  =0$ for $r\in\mathbb{Z}$

\item  If $d_{s}<j\leq d_{s-1}$ for some $1<s\leq p$ then $F_{r}Z^{j}\left(
\frak{g},\frak{g}\right)  =Z^{j}\left(  \frak{g},\frak{g}\right)  $ for $r\leq
q$, where%
\[
q=-\left[  pd_{p}+\left(  p-1\right)  \left(  d_{p-1}-d_{p}\right)
+..+s\left(  d_{s}-d_{s+1}\right)  +\left(  s-1\right)  \left(  j-d_{s}%
-1\right)  \right]
\]
\end{enumerate}
\end{lemma}

\begin{corollary}
Let $\frak{g}$ be an n-dimensional filiform Lie algebra. For $2r\leq\left(
j-1\right)  \left(  j-2p-2\right)  $ with $1\leq j\leq n-1$ we have
\[
F_{r}Z^{j}\left(  \frak{g},\frak{g}\right)  =Z^{j}\left(  \frak{g}%
,\frak{g}\right)
\]
\end{corollary}

Further, it can be proven ( see [96] or [54] ) that if $r\leq p-pj$, then
$F_{r}Z^{j}\left(  \frak{g},\frak{g}\right)  =Z^{j}\left(  \frak{g}%
,\frak{g}\right)  $. As a consequence, any derivation of the Lie algebra
$\frak{g}$ will map the space $F_{r}\frak{g}$ \ on $T_{r}\frak{g}$ for any
$r$. This leads to the equality given by Vergne, namely that for $r\leq-p$,
where $p$ is the nilindex of the algebra, we have $F_{r}H^{2}\left(
\frak{g},\frak{g}\right)  =H^{2}\left(  \frak{g},\frak{g}\right)  $. This
equality has been of importance in the study of the irreducible components of
the variety of filiform laws.\newline Now let $\frak{g}=L_{n}$ be the mopdel
filiform Lie algebra introduced in section 1. For this algebra, it is not
difficult to prove that its Lie algebra of derivations is $\left(
2n+1\right)  $-dimensional, where $\dim L_{n}=n+1$. Thus the dimension of the
cohomology space $H^{1}\left(  L_{n},L_{n}\right)  $ is also $n+1$, and from
this $\dim\;B^{2}\left(  L_{n},L_{n}\right)  =n^{2}$. The description of the
spaces $F_{0}Z^{2}\left(  L_{n},L_{n}\right)  $ is the key to construct its
characteristically nilpotent deformations. Let $\left\{  X_{0},..,X_{n}%
\right\}  $ be a basis of $L_{n}$ and define the cochians $\phi\left(
X_{0},X_{i}\right)  =X_{j}$ for 1$\leq i,j\leq n$. As they are cocycles, the
determination of the space  $Z^{2}\left(  L_{n},L_{n}\right)  $ is reduced to
the study of those cocycles which satisfy $\phi\left(  X_{0},X_{i}\right)  =0$
and preserve the natural graduation. In [54], the author construct the
following cocycles,  :%
\[
\psi_{k,s}\left(  X_{i},X_{i+1}\right)  =\left\{
\begin{array}
[c]{c}%
X_{s}\;if\;i=k\\
0\;if\;i\neq k
\end{array}
\right.
\]
Therefore, the remaining images are given by the relation
\[
\psi_{k,s}\left(  X_{i},X_{j}\right)  =\left(  -1\right)  ^{k-i}%
C_{k-i}^{j-k-1}\left(  adX_{0}\right)  ^{i+j-1-2k}X_{s}%
\]
Now these and the preceding cocycles describe the cohomology space $F_{0}%
Z^{2}\left(  L_{n},L_{n}\right)  $ completely :

\begin{proposition}
The cocycles $\phi_{i,j}$ and $\psi_{k,s}$ $\left(  i<j,\;s\leq2k+1\right)  $
form a basis of   $F_{0}Z^{2}\left(  L_{n},L_{n}\right)  $. 
\end{proposition}

\begin{corollary}
We have%
\[
\dim\;F_{0}H^{2}\left(  L_{n},L_{n}\right)  =\left\{
\begin{array}
[c]{c}%
\frac{3n^{2}-4n+1}{4}\;\text{for\ }n\equiv1\left(  \operatorname{mod}%
\,2\right)  \\
\frac{n^{2}-2n-4}{4}\;\text{for\ }n\equiv0\left(  \operatorname{mod}%
\,2\right)
\end{array}
\right.
\]
Moreover, a basis is given by the cohomology classes of $\psi_{k,s}$ for
$1\leq k\leq n,4\leq s\leq n$ whenever $s\geq2k+1$. 
\end{corollary}

Using the Chevalley cohomology of the Lie algebra $\frak{g}_{n}$ it can be
shown that the elements of the space $Z^{2}\left(  \frak{g}_{n},\frak{g}%
_{n}\right)  $ correspond to infinitesimal deformations of the algebra
$\frak{g}_{n}=\left(  \mathbb{C}^{n},\mu_{n}\right)  $ ( see [22], [25]). Let
$\psi$ be a cocycle and define the operation
\[
\left[  x,y\right]  _{\psi}:=\left[  x,y\right]  +\psi\left(  x,y\right)
,\;\;\;x,y\in\mathbb{C}^{n}%
\]
Then the deformation is linearly expandable if the previous operation
satisfies the Jacobi condition, i.e, defines a Lie algebra structure on
$\mathbb{C}^{n}$. Let $\psi
\in\bigoplus H_{i}^{2}\left(  L_{n},L_{n}\right)  =F_{1}H^{2}\left(
L_{n},L_{n}\right)  .$ Then the cocycle admits a decomposition $\psi
=\sum_{i=1}^{r}\psi_{i}$ with $\psi_{i}\in H_{i}^{2}\left(  L_{n}%
,L_{n}\right)  .$ The last nonzero component of this decomposition is called
the sill cocycle of $\psi$.\newline The idea used in [54] is to decompose the preceding basis into layers, where a layer $k_{0}$ contains those cocycles $\psi_{k,s}$ whose entry $k$ is $k_{0}$. Now a cocycles $\psi=\sum a_{k,s}\psi_{k,s}\in F_{0}H^{2}\left(L_{n},L_{n}\right)$ is called degenerate in the layer $k_{0}$ if all $a_{k_{0},s}$ are zero. If it is nondegenerate in this layer we choose $\psi_{k_{0},s_{0}}$ with $a-{k_{0},s_{0}}\neq 0$ of least class. This has been called the nondegeneracy class of $\psi$. Moreover, under the asumption that this last class is $r$, the layer $k_{0}$ is called special if $2k_{0}+r+1<s$ for any nonzero $a_{k,s}$ for which $k>k-{0}$. 

\begin{definition}
A nonzero cocycle $\psi\in F_{1}H^{2}\left(L_{n},L_{n}\right)$ is called regular if it is linearly expandable and satisfies one of the following conditions :
\begin{enumerate}
\item There exist two special layers in which the cocycle is nondegenerate with distincts nondegeneracy class.
\item The cocycle belongs to $F_{2}H^{2}\left(L_{n},L_{n}\right)$ and there exists a special layer $k_{0}$ of class $r$ such that $a_{k_{0},r+2+2k_{0}}\neq 0$ with $2k_{0}+r+2<s$ for those $a_{k,s}$ with $k>k_{0}$.
\end{enumerate}
\end{definition}

Provided with these cocycles, Khakimdjanov shows then the following

\begin{theorem}
Let $\psi$ be a regular cocycle in $F_{1}H^{2}\left(L_{n},L_{n}\right)$. Then the deformation $\left(L_{n}\right)_{\psi}$ is a CNLA.
\end{theorem}

\begin{corollary}
Let $S$ be the set of pairs $\left(  k,s\right)  $ \ of positive integers such
that $\frac{\left(  n-5\right)  }{2}+2k+1\leq s\leq n$ and $\psi=\sum_{\left(
k,s\right)  S}a_{k,s}\psi_{k,s}$. Let $s_{0}$ be the least integer such that
$s_{0}\geq\frac{n+1}{2}$. If one of the following conditions 

\begin{enumerate}
\item $n>8$ and $a_{1,s_{0}},a_{1,s_{0}+1}\neq0,$

\item $n\geq6,\;a_{1,s_{0}}=0$ and $a_{1,s_{0}+1},a_{1,s_{0}+2}\neq
0,$\newline holds, then $\left(  L_{n}\right)  _{\psi}$ is a CNLA.
\end{enumerate}
\end{corollary}

This and other corollaries contained in [54] allow to cosntruct large
families of CNLAs. The idea is to consider subsets of the basis given above
such that the elements of the linear envelope of this set gives lineraly
expandable cocycles. Imposing additional conditions on the coefficients, the
cocycles are made regular. It is remarked that there exist characteristically
nilpotent deformations of $L_{n}$ based on nonregular cocycles [54].
Moreover, the closure of the orbit corresponding to the set of CNLA of the
preceding corollary is a closed irreducible set of the variety $\frak{N}%
^{n+1}$ containing a nonempty Zariski-open subset formed by CNLAs. \newline Other results of the same nature due to this author are the follwoing :

\begin{lemma}
Let $\psi\in F_{1}H^{2}\left(  L_{n},L_{n}\right)  $ be a linearly expandable
nonzero cocycle. Then its sill cocycle $\psi_{r}$ is also linearly expandable.
\end{lemma}

Now let $\left(  L_{n}\right)  _{\psi}$ be a deformation with
$\psi\in F_{1}H^{2}\left(  L_{n},L_{n}\right)  .$ Let $\psi_{r}$ be the sill
cocycle of $\psi.$ Then the Lie algebra $\left(  L_{n}\right)  _{\psi_{r}}$ is
called the sill algebra of $\left(  L_{n}\right)  _{\psi}.$ The relation
between these two algebras is the crucial point to construct
characteristically nilpotent Lie algebras$\footnote{It is evident that the
infinitesimal deformations are filiform, for we have seen that the
characteristic sequence of the deformation is greater or equal than $c\left(
L_{n}\right)  ,$ and this is the maximal one.}$

\begin{theorem}
Let $\psi\in F_{1}H^{2}\left(  L_{n},L_{n}\right)  $ be a nonzero linearly
expandable cocycle. Then the Lie algebra $\left(  L_{n}\right)  _{\psi}$ is
characteristically nilpotent if and only if it is not isomorphic to its sill
algebra $\left(  L_{n}\right)  _{\psi_{r}}.$
\end{theorem}

From the theorem we obtain for example the following characteristic nilpotent
Lie algebras with basis $\left\{  X_{0},..,X_{2m}\right\}  $ and law
\[%
\begin{tabular}
[c]{cl}%
$\lbrack X_{0},X_{i}]=X_{i+1},$ & $i=1,..,2m-1$\\
$\lbrack X_{1},X_{i}]=X_{i+3},$ & $i=2,..,2m-3$\\
$\lbrack X_{i},X_{2m-i-1}]=\left(  -1\right)  ^{i+1}X_{2m}$ & $i=1,..,m-1$%
\end{tabular}
\]

For the nonfiliform Lie algebras the determination of characteristically
nilpotent Lie algebras is not so well structured. In fact, for any lower
characteristic sequence there will appear more naturally graded models than it
was the case in the filiform algebras. This construction allowed to obtain certain results on the structure
of the neighborhhods of filiform Lie algebras on the variety $\frak{N}^{n}$
[23], so it is of interest for the determination of the irreducible components
of the variety of filiform Lie algebra laws, thus for the variety
$\frak{N}^{n}$ itself. We mantain the notation for the cohomology introduced earlier.

\begin{lemma}
Let $s>r,s\neq2r.$ If there is a nonzero cocycle $\psi\in H_{s}^{2}\left(
L_{n},L_{n}\right)  $ belonging to $H_{s}^{2}\left(  L_{n},L_{n}\right)  \cap
B^{2}\left(  \left(  L_{n}\right)  _{\psi},\left(  L_{n}\right)  _{\psi
}\right)  ,$then this cocycle is unique ( up to multiples).
\end{lemma}

The proof is based on the structure of the algebra of derivations of a sill
algebra and is omitted here. It can be found in [23] and [26]. Now let
$A=\left(  L_{n}\right)  _{\psi}$ be a filiform algebra, where $\psi\in
Z^{2}\left(  L_{n},L_{n}\right)  \cap F_{1}H^{2}\left(  L_{n},L_{n}\right)  $
and $\psi_{r}$ denotes the sill cocycle of $\psi.$

\begin{lemma}
Let $n\geq8$ and $V$ an open set of $\frak{N}^{n}$ containing $A.$ Then there
exists a characteristically nilpotent Lie algebra in $V.$
\end{lemma}

Then we obtain immediately the following 
\begin{corollary}
For $\frak{N}^{n}\;\left(  n\geq7\right)  $ \ there exists an open set whose
elements are characteristically nilpotent Lie algebras.
\endproof
\end{corollary}

\section{Lie algebras of type $Q$ and its deformations}
In this section we use the other naturally graded filiform Lie algebra, $Q_{n}$, to obtain characteristically nilpotent Lie algebras in any dimension $n\geq 9$ and mixed characteristic sequence. This approach is perhaps not so natural, but it is based on an important property of "noncommutativity", which allows to obtain "easier" deformations. Combined with central extensions of special kind, we obtain the desired characteristically nilpotent deformations.
Let us concentrate on the Lie algebra $Q_{n}$. In contrast to $L_{n}$, it can only appear in even dimension. Thus the algebra $Q_{n}$ posesses a structural obstruction that forces its even-dimensionality. This obstruction is strongly related with the properties of the descending central sequence ${C^{k}Q_{n}}$.
\newline Let ${\omega_{1},..,\omega_{2m}}$ be the dual basis of the basis ${X_{1},..,X_{2m}}$ of $Q_{n}$. Then the Cartan-Maurer equations of this algebra are :

\begin{eqnarray*}
d\omega _{1} &=&d\omega _{2}=0 \\
d\omega _{j} &=&\omega _{1}\wedge \omega _{j-1},\;3\leq j\leq 2m-1 \\
d\omega _{2m} &=&\omega _{1}\wedge \omega _{2m-1}+\sum_{j=2}^{m}\left(
-1\right) ^{j}\omega _{j}\wedge \omega _{2m+1-j}
\end{eqnarray*}

In particular, the nonzero exterior product $\omega_{m}\wedge\omega_{m+1}$ shows that the ideal $C^{p-1}Q_{n}$, where $\left[ \frac{2m-1}{2}\right]$ and $n=2m-1$ of the central descending sequence is not abelian, while $C^{p}Q_{n}$ is abelian. This can be interpreted in the following manner: while $L_{n}$ has abelian commutator algebra $C^{1}L_{n}$, the model $Q_{n}$ is as far as possible from being an abelian algebra. This fact is important for deformation theory, as it can be interpreted in the sense that deforming $Q_{n}$ will be easier than deforming $L_{n}$.
\newline The previous property can be expressed in terms of centralizers :

\begin{eqnarray*}
C_{Q_{n}}\left( C^{p}Q_{n}\right)  &\supset &C^{p}Q_{n} \\
C_{Q_{n}}\left( C^{q}Q_{n}\right)  &\varsupsetneq &C^{q}Q_{n}
\end{eqnarray*}
for $n=2m-1$, $p=\left[ \frac{2m-1}{2}\right]$ and $1\leq q\leq p-1$.\newline We will say that $Q_{n}$ satisfies the centralizer property.\newline It is rather convenient to generalize this property to any naturally graded Lie algebra :

\begin{definition}
Let $\frak{g}_{n}$ be an $n$-dimensional, naturally graded nilpotent Lie algebra of nilindex $p$. Then $\frak{g}_{n}$ is called of type $Q$ if 
\begin{eqnarray*}
C_{\frak{g}_{n}}\left( C^{p}\frak{g}_{n}\right)  &\supset &C^{p}\frak{g}_{n} \\
C_{\frak{g}_{n}}\left( C^{q}\frak{g}_{n}\right)  &\varsupsetneq &C^{q}\frak{g}_{n}
\end{eqnarray*}
for $n=2m-1$, $p=\left[ \frac{p}{2}\right]$ and $1\leq q\leq p-1$.
\end{definition}

We are principally concerned with the Lie algebras of type $Q$ that are central extensions of the filiform Lie algebra $Q_{n}$, as well as other extensions.\newline Observe however that the index fixed in the previous definition is maximal, i.e, there do not exist Lie algebras which are "less abelian" with respect to the previous definition. The index, will depends only on the nilindex of the algebra, is very important and appears in apparently different contexts, such as the parabolic subalgebras [58].

\begin{theorem}
Let $\frak{n}$ be the nilradical of a standard Borel subalgebra $\frak{b}%
\left(  \Delta\right)  $ of a complex simple Lie algebra distinct from $G_{2}%
$. Then $\frak{n}$ satisfies the centralizer property.
\end{theorem}

The proof is an immediate consequence of the following result :

\begin{proposition}
Let $\frak{n}$ be the nilradical of a standard Borel subalgebra $\frak{b}%
\left(  \Delta\right)  $ of a complex simple Lie algebra distinct from $G_{2}%
$. Let $p=ht\left(  \delta\right)  $ be the height of the maximal root. Then
there exist roots $\alpha,\beta$  whose height is $[\frac{ht\left(  \delta\right)  }{2}]$ such that $\alpha+\beta$ is a positive root.
\end{proposition}

Thus we see that the classical theory provides a lot of naturally graded Lie algebras satisfying the centralizer property. However, it is usually unconvenient to manipulate these algebras, because of the great difference between its dimension and nilpotence class : the first is too high in comparison with the last.
\newline From the definition it follows also that a central extension $\frak{e}$ of $Q_{n}$ by $\mathbb{C}$ of type $Q$ cannot be filiform. This implies that the cocycle $\varphi\in H^{2}\left(Q_{n},\mathbb{C}\right)$ that defines the extension cannot be affine [17]. As a central extension of a filiform Lie algebra is filiform if and only if the cohomology class of $\varphi$ is affine, we conclude that for our special case, the extension $\frak{e}$ cannot be given by an affine cocycle.\newline
Let $\frak{e}\in E_{c,1}\left(Q_{n}\right)$ be an extension of type $Q$. As the nilindex is preserved, we conclude that the characteristic sequence of $\frak{e}$ must be lower than $\left(2m,1\right)$. Thus these algebras will play, in the set of Lie algebras with this characteristic sequence, the same role that $Q_{n}$ plays for the filiform algebras.\newline

Let $\overset{\sim }{E}_{c,1}\left( Q_{n}\right) =\left\{ \frak{e}\in
E_{c,1}\left( Q_{n}\right) \;|\;\frak{e}\text{ \ is of type }Q\right\}$. If $\frak{e}$ is any such element expresed over the basis ${X_{1},..,X_{2m+1}}$, it follows immediately from the definition of type $Q$ that $\frak{e}$ is naturally graded. The first $2m$ vectors are fixed in the natural graduation of the extension, thus $\frak{e}$ is completely determined once we know the position of the vector $X_{2m+1}$ in the graduation. The next lemma establishes that the positions are not arbitrary.

\begin{lemma}
Let $\frak{e}\in E_{c,1}\left( Q_{n}\right) $ be an extension. If $%
X_{2m+1}\in \frak{e}_{2t}$ $\left( 1\leq t\leq \left[ \frac{2m-1}{2}\right]
\right) $ then $\frak{e}$ is not naturally graded. In particular, $\frak{e}%
\notin \overset{\sim }{E}_{c,1}\left( Q_{n}\right)$.
\end{lemma}

It follows that the position of the vector $X_{2m+1}$ is only admissible if the graduation block is odd indexed. As we are not interested in split Lie algebras, we convene that $X_{2m+1}\notin\frak{e}_{1}$. Moreover, we define the depth $h$ of $X_{2m+1}$ like follows:
\[
h\left(X_{2m+1}\right)=t\; \text{if}\; X_{2m+1}\in\frak{e}_{2t+1},\; 1\leq t\leq \left[ \frac{2m-1}{2}\right]-1
\]
For convenience Lie algebras will be written usually in their contragradient representation. This will be of importance for the deformations, as linearly expandable cocycles are easier recognized when using this notation. Let ${\omega_{1},..,\omega_{2m+1}}$ be the dual basis to ${X_{1},..,X_{2m+1}}$ for the extension $\frak{e}\in E_{c,1}\left(Q_{n}\right)$. Then its Cartan-Maurer equations are :
\begin{eqnarray*}
d\omega _{1} &=&d\omega _{2}=0 \\
d\omega _{j} &=&\omega _{1}\wedge \omega _{j-1},\;3\leq j\leq 2m-1 \\
d\omega _{2m} &=&\omega _{1}\wedge \omega _{2m-1}+\sum_{j=2}^{m}\left(
-1\right) ^{j}\omega _{j}\wedge \omega _{2m+1-j}\\
d\omega_{2m+1} &=&\sum_{i,j}a^{ij}\omega _{i}\wedge \omega _{j},\;a^{ij}\in \mathbb{C}%
,\;i,j\geq 2
\end{eqnarray*}
where $d^{2}\omega_{2m+1}=0$. Then the determination of the extensions of type $Q$ of $Q_{n}$ reduces to the determination of the possible differential forms $d\omega_{2m+1}$. As known, the coefficient $a^{i,j}$ is given by a linear form over $\bigwedge^{2} Q_{n}$ which annihilates over $\Omega$.\newline

Let $\varphi _{ij}\in Hom\left( \bigwedge^{2}Q_{n},\mathbb{C}\right)
,\,2\leq i,j\leq 2m$, be defined by 
\begin{equation*}
\varphi _{ij}\left( X_{k},X_{l}\right) =\delta _{ik}\delta _{kl},\;\;\left(
X_{k},X_{l}\right) \in \frak{g}^{2}
\end{equation*}

\begin{lemma}
For $m\geq 4$ and $1\leq t\leq m-2$ \ the cochain  $\varphi
_{t}=\sum_{j=2}^{t+1}\left( -1\right) ^{j}\varphi _{j,3+2t-1}$ defines a
cocycle of $H^{2}\left( Q_{n},\mathbb{C}\right) $. If $\frak{g}_{\left(
m,t\right) }$ denotes the extension defined by $\varphi _{t}$, then $\frak{g}%
_{\left( m,t\right) }\in \overset{\sim }{E}_{c,1}\left( Q_{n}\right) $.
\end{lemma}
In particular, it follows from the proof [22] that the Cartan-Maurer equations of such an extension are 
\begin{eqnarray*}
d\omega _{1} &=&d\omega _{2}=0 \\
d\omega _{j} &=&\omega _{1}\wedge \omega _{j-1},\;3\leq j\leq 2m-1 \\
d\omega _{2m} &=&\omega _{1}\wedge \omega _{2m-1}+\sum_{j=2}^{m}\left(
-1\right) ^{j}\omega _{j}\wedge \omega _{2m+1-j}\\
d\omega_{2m+1} &=&\sum_{j=2}^{t+1}\left(-1\right)^{j}\omega _{j}\wedge \omega _{3-j+2t},\;1\leq t\leq m-2
\end{eqnarray*}

The family of extensions ( which is locally finite and depends on $m$ ) is
proven to be the class of algebras we are interested in, as follows from the
next 

\begin{proposition}
An extension $\frak{e}\in E_{c,1}\left( Q_{n}\right) $ is of type $Q$ if and
only if there exists a $t\in \left\{ 1,..,m-2\right\} $ such that $\frak{e}%
\simeq \frak{g}_{\left( m,t\right) }$.
\end{proposition}

Let $\widehat{H}^{2}\left( Q_{n},\mathbb{C}\right) =\left\{ \varphi \in
H^{2}\left( Q_{n},\mathbb{C}\right) \;|\;\frak{e}_{\varphi }\text{ is of
type }Q\right\} $, where $\frak{e}_{\varphi }$ is the extension defined by $%
\varphi $. The above result proves that $\dim \;\widehat{H}^{2}\left( Q_{n},%
\mathbb{C}\right) =m-1$, where $n=2m-1$. Moreover, the type of the extension $\frak{g}_{\left( m,t\right) }$ satisfies
\begin{eqnarray*}
p_{1} &=&p_{2t+1}=2 \\
p_{j} &=&1\text{ \ if }1\leq j\leq 2m-1,j\notin \left\{ 1,2m+1\right\} 
\end{eqnarray*}

As we have seen, the structure of the extensions $\frak{g}_{\left( m,t\right) }$ is very similar, in the sense that the differential form $d\omega_{2m+1}$ has a precise form which depends only on the depth of the ( added ) vector $X_{2m+1}$ dual to $X_{2m+1}$. 
\newline Now a construction method for characteristically nilpotent Lie algebras is given.  These deformations will be also interpretable in term of the graded cohomology spaces $H_{k}^{2}\left(  \frak{g},\frak{g}\right)$ associated to the lie algebra $\frak{g}$. The results given here are more widely covered in [22] :  

\begin{notation}
Let $\frak{g}$ be a $n$-dimensional Lie algebra defined over the basis 
 $\left\{  X_{1},..,X_{n}\right\}  $ and let  $Der\left(  \frak{g}\right)  $
be its algebra of derivations. If $f\in Der\left(  \frak{g}\right)  $,
we will use the notation
\[
f\left(  X_{i}\right)  =\sum_{j=1}^{n}f_{i}^{j}X_{j},\;1\leq i\leq n
\]
\end{notation}

We consider the following cocycle ( class ) for the Lie algebras $\frak{g}_{\left(m,t\right)}$ and $t\geq 2$ : 
\[
\varphi_{m,t}\left(  X_{2},X_{3+j}\right)  =X_{2t+2+j},\;0\leq j\leq2m-2t-2
\]
The reason for excluding the value $t=1$ lies in the simplicity of its last differential form. For these algebras special cocycles have to be considered [22] :

\begin{lemma}
For  $\ m\geq5,\;1\leq t\leq m-2$ \ $\varphi_{m,t}$ is linearly expandable.
\end{lemma}

\begin{proposition}
For  $m\geq5,\;1\leq t\leq m-2$ the Lie algebra $\frak{g}_{\left(
m,t\right)  }+\varphi_{m,t}$ is characteristically nilpotent.
\end{proposition}

Note that the cocycle which defines the deformation $\frak{g}_{\left(
m,t\right)  }+\varphi_{m,t}$ is chosen such that the incorporated 
brackets do not change the exterior differential of the system.
The cocycle $\varphi_{m,t}$ admits the following cohomological interpretation :

\begin{proposition}
For $t\geq3$ \ let $\psi\in H_{2t-2}^{2}\left(  \frak{g}_{\left(  m,t\right)
},\frak{g}_{\left(  m,t\right)  }\right)  \;$ be a cocycle that satisfies

\begin{enumerate}
\item $\forall\;X\in Z\left(  \frak{g}_{\left(  m,t\right) }\right)  $
such that  $h\left(  X\right)  =t$, we have  $\psi\left(  X,\frak{g}_{\left(
m,t\right)  }\right)  =\{0\}$ and $X\notin im\left(  \psi\right)  $

\item  If $X\in\frak{g}_{\left(  m,t\right)  }$ is such that there exists an $Y\in
Z\left(  \frak{g}_{\left(  m,t\right)  }\right)  $ with $h\left(  Y\right)
=t$ and $Y\notin im\,\,ad\left(  X\right)  $, then $\psi\left(
X,C^{1}\frak{g}_{\left(  m,t\right)  }\right)  =\{0\}.$\newline Then
\[
\psi=\sum_{\substack{2\leq i\leq t+2\\3\leq j\leq2m-3}}\lambda_{ij}\psi
_{ij}\;\;\left(  \lambda_{ij}\in\mathbb{C}\right)
\]
where
\[
\psi_{ij}\left(  X_{i},X_{j}\right)  =X_{i+j+1},\;\;i+j\leq2m\text{ \ \ }%
\]
\end{enumerate}
\end{proposition}

Writing 
\[
\widehat{H}_{2t-2}^{2}\left(  \frak{g}_{\left(  m,t\right)  }%
,\frak{g}_{\left(  m,t\right)  }\right)  =\left\{  \psi\;|\;\psi\;\text{
satisfies\ }1)\text{ and\ }2)\right\}
\]
we isolate the cohomology classes that give the desired deformations :

\begin{corollary}
A cocycle $\psi\in$ $\widehat{H}_{2t-2}^{2}\left(  \frak{g}_{\left(
m,t\right)  },\frak{g}_{\left(  m,t\right)  }\right)  $ such that
$\psi\left(  C^{1}\frak{g}_{\left(  m,t\right)  },C^{1}\frak{g}_{\left(
m,t\right)  }\right)  =\{0\}$ is linearly expandable if and only if
$\psi=\lambda\varphi_{m,t}\;\left(  \lambda\in\mathbb{C}\right)
$.\ Moreover, $\frak{g}_{\left(  m,t\right)  }+\lambda\varphi
_{m,t}\simeq\frak{g}_{\left(  m,t\right)  }+\varphi_{m,t}$ for any
$\lambda\neq0$. 
\end{corollary}

From the corollary we deduce that $\varphi_{m,t}$ is fixed, up to
multiples, by the restriction property to the derived subalgebra.

\begin{theorem}
Let $\psi\in\widehat{H}_{2t-2}^{2}\left(  \frak{g}_{\left(  m,t\right)  }%
,\frak{g}_{\left(  m,t\right)  }\right)  $ be a linearly 
expandale cocycle. Then the algebra $\frak{g}_{\left(  m,t\right)  }%
+\psi$ is characteristically nilpotent.
\end{theorem}

Any supplementary deformation to the one defined by the cocycle
$\varphi_{m,t}$ changes the law  $\frak{g}_{\left(m,t\right)}$ in the same
way as  $\varphi_{m,t}$, so that it does not alter the conditions on the derivations.
Further, we determine certain central extensions of the algebras $%
\frak{g}_{\left( m,t\right) }$ obtained before. Observe that the
characteristic of an extension of $\frak{g}_{\left( m,t\right) }$ by $%
\mathbb{C}$ can be either $\left( 2m-1,1,1,1\right) $ or $\left(
2m-1,2,1\right) $. The first one is not interesting for our purposes, as it
is linear, while the second one is mixed\footnote{%
A characteristic sequence $c\left( \frak{g}\right) $ is called mixed if
there are two or more Jordan blocks of dimension $\geq 2$.}. 

Let $\mathcal{G}_{2}^{1}=\left\{ \frak{g}_{\left( m,t\right) }\;|\;m\geq
4,\;1\leq t\leq m-2\right\} $. For any fixed $m$ and $t$ we define 
\begin{equation*}
E_{c,1}^{1}\left( \frak{g}_{\left( m,t\right) }\right) =\left\{ \frak{e}\in
E_{c,1}\left( \frak{g}_{\left( m,t\right) }\right) \;|\;\frak{e}\text{ is of
type }Q\text{ and }h\left( X_{2m+2}\right) =h\left( X_{2m+1}\right)
+1\right\} 
\end{equation*}
where $\left\{ X_{1},..,X_{2m+1}\right\} $ is a basis of $\frak{g}_{\left(
m,t\right) }$, $\left\{ X_{1},..,X_{2m+2}\right\} $ a basis of $\frak{e}$
and $h$ is the depth funtion.

\begin{proposition}
Let $t\geq 2$ and $\frak{g}_{\left( m,t\right) }\in \mathcal{G}_{2}^{1}$.
Then an extension $\frak{e}\in E_{c,1}\left( \frak{g}_{\left( m,t\right)
}\right) $ belongs to $E_{c,1}^{1}\left( \frak{g}_{\left( m,t\right)
}\right) $ if and only if its structural equations are
\begin{eqnarray*}
d\omega _{1} &=&d\omega _{2}=0 \\
d\omega _{j} &=&\omega _{1}\wedge \omega _{j-1},\;3\leq j\leq 2m-1 \\
d\omega _{2m} &=&\omega _{1}\wedge \omega _{2m+1}+\sum_{j=2}^{m}\left(
-1\right) ^{j}\omega _{j}\wedge \omega _{2m+1-j} \\
d\omega _{2m+1} &=&\sum_{j=2}^{t+1}\left( -1\right) ^{j}\omega _{j}\wedge
\omega _{3+2t-j} \\
d\omega _{2m+2} &=&\omega _{1}\wedge \omega _{2m+1}+\sum_{j=2}^{t+1}\left(
-1\right) ^{j}\left( t+2-j\right) \omega _{j}\wedge \omega _{4+2t-j}
\end{eqnarray*}
\end{proposition}

An extension $\frak{e}$ with the previous Cartan-Maurer equations will be
denoted by $\frak{g}_{\left( m,t\right) }^{1}$.
Observe that the case $t=1$ has been excluded from the proposition. The
reason is that, by the simplicity of the differential form $d\omega _{2m+1}$%
, in this case there are two possible extensions.

\begin{lemma}
For $m\geq 4$, $\frak{e}\in E_{c,1}\left( \frak{g}_{\left( m,1\right)
}\right) $ belongs to $E_{c,1}^{1}\left( \frak{g}_{\left( m,1\right)
}\right) $ if the structural equations of $\frak{e}$ over a basis $\left\{
\omega _{1},..,\omega _{2m+2}\right\} $ are 
\begin{eqnarray*}
d\omega _{1} &=&d\omega _{2}=0 \\
d\omega _{j} &=&\omega _{1}\wedge \omega _{j-1},\;3\leq j\leq 2m-1 \\
d\omega _{2m} &=&\omega _{1}\wedge \omega _{2m+1}+\sum_{j=2}^{m}\left(
-1\right) ^{j}\omega _{j}\wedge \omega _{2m+1-j} \\
d\omega _{2m+1} &=&\omega _{2}\wedge \omega _{3} \\
d\omega _{2m+2} &=&\omega _{1}\wedge \omega _{2m+1}+\omega _{2}\wedge \omega
_{4}+k\omega _{2}\wedge \omega _{2m+1},\;k=0,1
\end{eqnarray*}
\end{lemma}

The proof is analogous to the preceding one. The reason for the existence of
the second extension is the weakness of the restrictions imposed by the
differential form $d\omega _{2m+1}$. For higher depths the existence of
additional exterior products in the adjoined form $d\omega _{2m+2}$ is not
compatible with its closure $d^{2}\omega _{2m+2}=0$. 

\begin{notation}
For $k=0$ the extension is denoted by $\frak{g}_{\left( m,1\right) }^{1}$,
and for $k=1$ by $\frak{g}_{\left( m,1\right) }^{2}$.
\end{notation}

As known, the set of nilpotent Lie algebras $\frak{g}$ of a given dimension $%
n$ and characteristic sequence $c\left( \frak{g}\right) $ is denoted by $%
\delta U_{c\left( \frak{g}\right) }^{n}$ [3]. Now let $E_{c,2}\left(
Q_{n}\right) $ be the set of central extensions of $Q_{n}$ by $\mathbb{C}^{2}
$. The following result shows that we have obtained practically all the
extensions that interest us. \newline
Let $\frak{g}_{\left( m,0\right) }^{1+k}\;\left( k=0,1\right) $ be the Lie
algebras with structural equations
\begin{eqnarray*}
d\omega _{1} &=&d\omega _{2}=0 \\
d\omega _{j} &=&\omega _{1}\wedge \omega _{j-1},\;3\leq j\leq 2m-1 \\
d\omega _{2m} &=&\omega _{1}\wedge \omega _{2m+1}+\sum_{j=2}^{m}\left(
-1\right) ^{j}\omega _{j}\wedge \omega _{2m+1-j} \\
d\omega _{2m+1} &=&0 \\
d\omega _{2m+2} &=&\omega _{1}\wedge \omega _{2m+1}+k\omega _{2}\wedge
\omega _{2m+1}
\end{eqnarray*}

\begin{theorem}
For $n=2m-1,\;m\geq 4$ the following identity holds :
\begin{equation*}
E_{c,2}\left( Q_{n}\right) \cap \delta U_{\left( 2m-1,2,1\right)
}^{2m+2}=\bigcup_{j=2}^{m-2}\mathcal{O}\left( \frak{g}_{\left( m,t\right)
}^{1}\right) \cup \mathcal{O}\left( \frak{g}_{\left( m,1\right) }^{2}\right)
\cup \mathcal{O}\left( \frak{g}_{\left( m,0\right) }^{1+k}\right) ,\;k=0,1
\end{equation*}
where $\mathcal{O}\left( \frak{g}\right) $ denotes the orbit of $\frak{g}$
by the action of the general linear group.
\end{theorem}

Any extension of $Q_{n}$ by $\mathbb{C}^{2}$ must have characteristic
sequence $\left( 2m-1,1,1,1\right) $ or $\left( 2m-1,2,1\right) $ if it
preserves the nilindex. Observe however that for the first sequence, the
split algebra $Q_{n}\oplus \mathbb{C}$ cannot generate a nonsplit central
extension.
Now it is convenient to introduce some notation. For $1\leq
t\leq m-2$ we can write the algebras $\frak{g}_{\left( m,t\right) }^{1}$
formally as 
\begin{equation*}
\frak{g}_{\left( m,t\right) }^{1}=\frak{g}_{\left( m,t\right) }+d\overset{-}{%
\omega }_{m,t}
\end{equation*}
where 
\begin{equation*}
d\overset{-}{\omega }_{m,t}=\omega _{1}\wedge \omega
_{2m+1}+\sum_{j=2}^{t+1}\left( -1\right) ^{j}\left( t+2-j\right) \omega
_{j}\wedge \omega _{4-j+2t}
\end{equation*}
is called extensor of type I.

\subsection{Deformations of $\frak{g}_{\left(  m,t\right)  }^{1}\;\left(
t\geq2\right)  $}

Let  $\frak{g}_{\left(  m,t\right)  }^{1}$ and
consider an extensor of type I \ $d\overset{-}{\omega}_{m,t}$. We know that
$\frak{g}_{\left(  m,t\right)  }^{1}=\frak{g}_{\left(  m,t\right)  }%
+d\overset{-}{\omega}_{m,t}$.
\newline Consider a cocycle $\psi\in H^{2}\left(  \frak{g}_{\left(  m,t\right)
}^{1},\frak{g}_{\left(  m,t\right)  }^{1}\right)  $ defined by
\[
\psi\left(  X_{i},X_{j}\right)  =\left\{
\begin{array}
[c]{c}%
\varphi_{m,t}\left(  X_{i},X_{j}\right)  \;\text{if }1\leq i,j\leq2m+1\\
0\;\;\text{if\ \ }i=2m+2\text{ or }j=2m+2
\end{array}
\right.
\]
$\psi$ is clearly a prolongation by zeros of the cocycle $\varphi
_{m,t}$; it will be convenient to preserve the notation
$\varphi_{m,t}$ to denote $\psi$, whenever there is no ambiguity.
In the previous section we saw that the adjoined extensors have no influence on 
the characetristic nilpotence of the deformation $\frak{g}_{\left(
m,2\right)  }+\varphi_{m,2}$. This property is in fact generalizable to any $t\geq3$ :

\begin{proposition}
For any $m\geq4,\;1\leq t\leq m-2$ the cocycle $\varphi_{m,t}\in H^{2}\left(
\frak{g}_{\left(  m,t\right)  }^{1},\frak{g}_{\left(  m,t\right)  }%
^{1}\right)  $ is linearly expandable.
\end{proposition}

\begin{corollary}
For any $m\geq4,\;1\leq t\leq m-2$ the following identity holds%
\[
\left(  \frak{g}_{\left(  m,t\right)  }+d\overset{-}{\omega}_{m,t}\right)
+\varphi_{m,t}=\left(  \frak{g}_{\left(  m,t\right)  }+\varphi
_{m,t}\right)  +d\overset{-}{\omega}_{m,t}%
\]
\end{corollary}

\begin{theorem}
For $m\geq4,\;1\leq t\leq m-2$ the Lie algebra $\frak{g}_{\left(
m,t\right)  }^{1}+\varphi_{m,t}$ is characteristically nilpotent
\end{theorem}

The cocycles $\varphi_{m,t}$ are a
special case of a more wide family of cocycles of the subspace
$H_{2t-2}^{2}\left(  \frak{g}_{\left(  m,t\right)  }^{1},\frak{g}_{\left(
m,t\right)  }^{1}\right)  $ :

\begin{lemma}
If $\psi\in H_{2t-2}^{2}\left(  \frak{g}_{\left(  m,t\right)  }^{1}%
,\frak{g}_{\left(  m,t\right)  }^{1}\right)  $ is a prolongation by
zeros of a cocycle $\varphi\in H_{2t-2}^{2}\left(  \frak{g}_{\left(
m,t\right)  },\frak{g}_{\left(  m,t\right)  }\right)  $, then
$\psi$ satisfies the conditions

\begin{enumerate}
\item $\forall\;X\in Z\left(  \frak{g}_{\left(  m,t\right)  }^{1}\right)  $
such that $h\left(  X\right)  =\frac{2t+1}{2}$, we have $\psi\left(
X,\frak{g}_{\left(  m,t\right)  }^{1}\right)  =\{0\}$ \ and $X\notin im\left(
\psi\right)  .$

\item $\forall\;X\in Z^{2}\left(  \frak{g}_{\left(  m,t\right)  }%
^{1}\right)  $ such that $h\left(  X\right)  =t$, we have $\psi\left(
X,\frak{g}_{\left(  m,t\right)  }^{1}\right)  =\{0\}$ \ and \ $X\notin
im\left(  \psi\right)  .$

\item  If $X\in\frak{g}_{\left(  m,t\right)  }^{1}$ is such that there exists an
$Y\in Z^{2}\left(  \frak{g}_{\left(  m,t\right)  }^{1}\right)  $ with
$h\left(  Y\right)  =t$ and $Y\notin im\,ad\left(  X\right)  $, then
$\psi\left(  X,C^{1}\frak{g}_{\left(  m,t\right)  }^{2,1}\right)  =\{0\}$.
\end{enumerate}
\end{lemma}

\begin{proposition}
A cocycle $\psi\in H_{2t-2}^{2}\left(  \frak{g}_{\left(  m,t\right)  }%
^{1},\frak{g}_{\left(  m,t\right)  }^{1}\right)  $ is a prolongation
by zeros of a cocycle $\varphi\in H_{2t-2}^{2}\left(  \frak{g}_{\left(
m,t\right)  },\frak{g}_{\left(  m,t\right)  }\right)  $ if and only
if it satisfies conditions $1),2),3)$.$\blacksquare$
\end{proposition}

We note
\[
\widehat{H}_{2t-2}^{2}\left(  \frak{g}_{\left(  m,t\right)  }^{1}%
,\frak{g}_{\left(  m,t\right)  }^{1}\right)  =\left\{  \psi\;|\;\psi\text{
satisfies }1),2)\text{ and }3)\right\}
\]

\begin{corollary}
A cocycle $\psi\in H_{2t-2}^{2}\left(  \frak{g}_{\left(  m,t\right)  }%
^{1},\frak{g}_{\left(  m,t\right)  }^{1}\right)  $ is a prolongation
by zeros of $\varphi_{m,t}$ if and only if the restriction of $\psi$ to
the derived subalgebra $C^{1}\frak{g}_{\left(  m,t\right)  }^{1}$ is
identically zero.
\end{corollary}

\begin{theorem}
Let $\psi\in\widehat{H}_{2t-2}^{2}\left(  \frak{g}_{\left(  m,t\right)
}^{1},\frak{g}_{\left(  m,t\right)  }^{1}\right)  $ be linearly 
expandable. Then $\frak{g}_{\left(  m,t\right)  }%
^{1}+\psi$ is characteristically nilpotent.
\end{theorem}

These results can be resumed graphically. 
We introduce the following notations [22] :%
\begin{align*}
M_{m,1}^{1}\left(  \frak{g}_{\left(  m,t\right)  }+\varphi_{m,t}\right)    &
=\frak{g}_{\left(  m,t+1\right)  }+\varphi_{m,t+1}\\
D_{1,t}^{1}\left(  \frak{g}_{\left(  m,t\right)  }+\varphi_{m,t}\right)    &
=\frak{g}_{\left(  m+1,t\right)  }+\varphi_{m+1,t}\\
d\overline{\omega}_{m,t}\left(  \frak{g}_{\left(  m,t\right)  }+\varphi
_{m,t}\right)    & =\frak{g}_{\left(  m,t\right)  }^{1}+\varphi_{m,t}\\
M_{m,1}^{2}\left(  \frak{g}_{\left(  m,t\right)  }^{1}+\varphi_{m,t}\right)
& =\frak{g}_{\left(  m,t+1\right)  }^{1}+\varphi_{m,t+1}\\
D_{1,t}^{2}\left(  \frak{g}_{\left(  m,t\right)  }^{1}+\varphi_{m,t}\right)
& =\frak{g}_{\left(  m+1,t\right)  }^{1}+\varphi_{m+1,t}%
\end{align*}
for $m\geq4$ and $1\leq t\leq m-2$.

\begin{figure}[h]
\vspace{2cm}
\begin{center}
\setlength{\unitlength}{0.60cm}
\begin{picture}(13,9)
\thicklines
\put(0,0){\framebox(7,7){}}
\put(0,7){\line(1,1){3}}
\put(7,7){\line(1,1){3}}
\put(7,0){\line(1,1){3}}
\put(3,10){\line(1,0){7}}
\put(10,3){\line(0,1){7}}
\thinlines
\put(3,3){\framebox(7,7){}}
\put(0,0){\line(1,1){3}}
\put(3,0){\vector(1,0){2}}
\put(5,3){\vector(1,0){2}}
\put(0.5,0.5){\vector(1,1){2}}
\put(7.5,0.5){\vector(1,1){2}}
\put(3,7){\vector(1,0){2}}
\put(6,10){\vector(1,0){2}}
\put(0.5,7.5){\vector(1,1){2}}
\put(7.5,7.5){\vector(1,1){2}}
\put(0,2){\vector(0,1){2}}
\put(7,2){\vector(0,1){2}}
\put(3,5){\vector(0,1){2}}
\put(10,5){\vector(0,1){2}}
\put(-2,-0.5){\makebox(0,0){$\frak{g}_{\left(m,t\right)}+\varphi_{m,t}$}}
\put(9,-0.5){\makebox(0,0){$\frak{g}_{\left(m,t+1\right)}+\varphi_{m,t+1}$}}
\put(2.7,2.3){\makebox(0,0){$\frak{g}_{\left(m+1,t\right)}+\varphi_{m+1,t}$}}
\put(10.7,2.3){\makebox(0,0){$\frak{g}_{\left(m+1,t+1\right)}+\varphi_{m+1,t+1}$}}
\put(-2,7.1){\makebox(0,0){$\frak{g}_{\left(m,t\right)}^{1}+\varphi_{m,t}$}}
\put(9,7.1){\makebox(0,0){$\frak{g}_{\left(m,t+1\right)}^{1}+\varphi_{m,t+1}$}}
\put(2.3,11){\makebox(0,0){$\frak{g}_{\left(m+1,t\right)}^{1}+\varphi_{m+1,t}$}}
\put(10.7,11){\makebox(0,0){$\frak{g}_{\left(m+1,t+1\right)}^{1}+\varphi_{m+1,t+1}$}}
\put(4,0.5){\makebox(0,0){$M_{m,1}^{1}$}}
\put(8.5,1){\makebox(0,0){$D_{1,t+1}^{1}$}}
\put(1.5,1){\makebox(0,0){$D_{1,t}^{1}$}}
\put(5,3.5){\makebox(0,0){$M_{m+1,1}^{1}$}}
\put(1,4){\makebox(0,0){$d\overset{-}{\omega}_{m,t}$}}
\put(8,4){\makebox(0,0){$d\overset{-}{\omega}_{m,t+1}$}}
\put(4,6){\makebox(0,0){$d\overset{-}{\omega}_{m+1,t}$}}
\put(10.5,6){\makebox(0,0){$d\overset{-}{\omega}_{m+1,t+1}$}}
\put(1.5,8){\makebox(0,0){$D_{1,t}^{2}$}}
\put(4,7.5){\makebox(0,0){$M_{m,1}^{2}$}}
\put(6,10.5){\makebox(0,0){$M_{m+1,1}^{2}$}}
\put(8.5,8){\makebox(0,0){$D_{1,t+1}^{2}$}}

\end{picture}
\caption{\label{cube}}
\end{center}
\end{figure}

\begin{theorem}
For $m\geq 4$ and $1\leq t\leq m-2$ the faces of the following cube are commutative diagrams.
\end{theorem}

\section{Nilpotent Lie algebras and rigidity}

Let $\frak{L}^{n}$ be the algebraic variety of complex Lie algebra laws on
$\mathbb{C}^{n}$. Each open orbit of the natural action of $GL\left(
n,\mathbb{C}\right)  $ on $\frak{L}^{n}$ gives, considering its Zarisky closure, an
irreducible component of $\frak{L}^{n}$ .\ Therefore, only a finite number of those
orbits exists; or, equivalently, only a finite number of isomorphism classes
of Lie algebras with open orbit. The first results about rigid Lie algebras are due to Gerstenhaber [38],
Nijenhuis and Richardson [78].\ The last two authors have transformed the
topological problems related to rigidity into cohomological problems, proving
in that an algebra is rigid if the second group in the Chevalley cohomology is
trivial. This theorem allows the construction of examples or rigid Lie
algebras and is used in the proof that semisimple algebras are
rigid.\ However, the existence of rigid Lie algebras whose second cohomology
group is non trivial shows that the cohomological viewpoint is not fully
satisfactory in the study of rigidity. 

\begin{definition}
A Lie algebra $\frak{g}$ is called decomposable if it can be written
\[
\frak{g}=\frak{s}\oplus\frak{t}\oplus\frak{n}%
\]
where $\frak{s}$ is a Levi subalgebra, $\frak{n}$ the nilradical and
$\frak{t}$ an abelian subalgebra whose elements are $ad$-semisimple and which
satisfies $\left[  \frak{s},\frak{t}\right]  =0.$
\end{definition}

The abelian subalgebra $T$ of $Der\frak{g}$ defined by
\[
T=\{adX,\quad X\in\frak{t}\}
\]
is called, following Malcev, an exterior torus on $\frak{g}$ . It is called
maximal torus, if it is maximal for the inclusion. Malcev has proved that all
maximal torus are pairwise conjugated, thus they have the same dimension
called the \textit{rank} of $\frak{g}$ and noted $r(\frak{g}).$

\begin{theorem}
Rigid Lie algebras are algebraic
\end{theorem}

\subsection{Roots system associated to a rigid solvable Lie algebra}

Let $\mu_{0}$ be a solvable decomposable law on $\mathbb{C}^{n}$.\ We fix a
maximal exterior torus $T.$ Let $X$ be a non-zero vector such that
$ad_{\mu_{0}}X$ belongs to $T$ .

\begin{definition}
\ We say that $X$ is \textbf{regular} if the dimension of
\[
V_{0}\left(  X\right)  =\left\{  Y\quad\text{such that \quad}\mu_{0}\left(
X,Y\right)  =0\right\}
\]
is minimal ; that is, $\dim V_{0}\left(  X\right)  \leq\dim V_{0}\left(
Z\right)  $ for all $Z$ such that $ad_{\mu_{0}}Z$ belongs to $T$.
\end{definition}

\begin{definition}
Suppose that $\mu_{0}$ is not nilpotent. The \textbf{root system} of $\mu_{0}$
associated to $\left(  X,Y_{1},...,Y_{n-p},X_{1},...,_{p-1}\right)  $ is the
linear system (S) defined by the following equations :

$x_{i}+x_{j}=x_{k}$ if the $X_{k}$-coordinate of $\mu_{0}\left(  X_{i}%
,X_{j}\right)  $ is non-zero.

$y_{i}+y_{j}=y_{k}$ if the $Y_{k}$-coordinate of $\mu_{0}\left(  Y_{i}%
,Y_{j}\right)  $ is non-zero.

$x_{i}+y_{j}=y_{k}$ if the $Y_{k}$-coordinate of $\mu_{0}\left(  X_{i}%
,X_{j}\right)  $ is non-zero.

$y_{i}+y_{j}=x_{k}$ if the $X_{k}$-coordinate of $\mu_{0}\left(  Y_{i}%
,Y_{j}\right)  $ is non-zero.
\end{definition}

In these notations we state.

\begin{theorem}
If $rank\left(  S\right)  \neq\dim\left(  I_{0}\right)  -1$, the law $\mu_{0}$
is not rigid.
\end{theorem}

\begin{corollary}
If $\mu_{0}$ is rigid, the rank of a root system for $\mu_{0}$ is independent
of the basis $\left(  X,Y_{1},...,Y_{n-p},X_{1},...,X_{p-1}\right)  $ used for
its definition.
\end{corollary}

\begin{corollary}
If $\mu_{0}$ is rigid, there is regular vector $X$ such that $ad_{\mu_{0}}X $
is diagonal and its eigenvalues are integers.
\end{corollary}

Let $\frak{R}_{n}$ be the variety of $n$-dimensional solvable Lie algebras. The principal structure theorem referring to rigid Lie algebras was proven by Carles in [18] :

\begin{theorem}
Any Lie algebra $\frak{g}$ \ which is rigid in either $\frak{L}^{n}$ or
$\frak{R}^{n}$ is algebraic and belongs to one of the following cases 

\begin{enumerate}
\item  The radical $Rad\left(  \frak{g}\right)  $ is not nilpotent and $\dim
Der\left(  \frak{g}\right)  =\dim\frak{g}$ ( if moreover $co\dim C^{1}%
\frak{g}>1$, the algebra is complete )

\item  The radical is nilpotent and satisfies one of the following conditions

\begin{enumerate}
\item $\frak{g}$ is perfect;

\item $\frak{g}$ is the direct product of $\mathbb{C}$ by a rigid perfect Lie
algebra whose derivations are inner;

\item $\frak{g}$ is non-perfect, has no direct abelain factor and is of rank
zero; morover, for any ideal of codimension one is also of rank zero.
\end{enumerate}
\end{enumerate}
\end{theorem}

\begin{corollary}
Any Lie algebra $\frak{g}$ rigid in $\frak{R}^{n}$ is algebraic and satisfies
one of the following conditions

\begin{enumerate}
\item $\dim Der\left(  \frak{g}\right)  =\dim\frak{g\;}$( if moreover $co\dim
C^{1}\frak{g}>1$, the algebra is complete )

\item $\frak{g}$ is characteristically nilpotent, as well as any of its
codimension one ideals.
\end{enumerate}
\end{corollary}

From the structure of the derivations for filiform Lie algebras, as found for example in [40], it follows easily that none filiform Lie algebra can be rigid in $\frak{L}^{n}$ or $\frak{R}^{n}$; by Carles' theorem, such an algebra would be characteristically nilpotent, and a contradiction with the dimension formulas is served. Thus the counterexamples, if any, must be searched within the nonfiliform Lie algebras. This would give an effective answer To Vergne's conjeture ( 1970 ) 

\begin{conjecture}
For any $n\neq1$ there do not exist nilpotent Lie algebras which are rigid in
$\frak{L}^{n}$ or $\frak{R}^{n}$.
\end{conjecture}

Recently we have found another curious relation between CNLA and rigid algebras. If we consider the Lie algebra $\frak{g}_{\left(  m,m-1\right)  }$ $\left(  m\geq3\right)  $
defined by the equations
\begin{align*}
d\omega_{1} &  =d\omega_{2}=0\\
d\omega_{j} &  =\omega_{1}\wedge\omega_{j-1},\;3\leq j<2m\\
d\omega_{2m+1} &  =\sum_{j=2}^{\left[  \frac{2m+1}{2}\right]  }\left(
-1\right)  ^{j}\;\omega_{j}\wedge\omega_{2m+1-j}%
\end{align*}
it is immediate to see that its characteristic sequence is $\left(  2m-1,1,1\right)$ and
its rank is two. Then there exists deformations which are isomorphic to the nilradical of a solvable rigid law, as gives the 

\begin{proposition}
The solvable Lie algebras \ $\frak{r}_{\left(  m,m-1\right)  }$ \ $\left(
m\geq3\right)  $ defined by the equations%
\begin{align*}
d\omega_{1} &  =\omega_{2m+2}\wedge\omega_{1}\\
d\omega_{2} &  =\left(  2m-3\right)  \omega_{2m+2}\wedge\omega_{2}\\
d\omega_{j} &  =\omega_{1}\wedge\omega_{j-1}+\left(  2m-5+j\right)
\omega_{2m+2}\wedge\omega_{j},\;3\leq j\leq2m-1\\
d\omega_{2m} &  =\omega_{1}\wedge\omega_{2m-1}+\omega_{2}\wedge\omega
_{3}+\left(  4m-5\right)  \omega_{2m+2}\wedge\omega_{2m}\\
d\omega_{2m+1} &  =\sum_{j=2}^{\left[  \frac{2m+1}{2}\right]  }\left(
-1\right)  ^{j}\;\omega_{j}\wedge\omega_{2m+1-j}+\left(  6m-9\right)
\omega_{2m+2}\wedge\omega_{2m+1}\\
d\omega_{2m+2} &  =0
\end{align*}
are rigid and complete. Moreover, their nilradical has codimension one and is
isomorphic to the Lie algebra $\frak{g}_{\left(  m,m-1\right)  }+\psi,$ where
$\psi\in H^{2}\left(  \frak{g}_{\left(  m,m-1\right)  },\frak{g}_{\left(
m,m-1\right)  }\right)  $ is the linearly expandable cocycle defined by
$\psi\left(  X_{2},X_{3}\right)  =X_{2m}.$
\end{proposition}

These algebras are a particular case of rigid Lie algebras whose nilradical
has codimension one, characteristic sequence $\left(  2m-1,1,1\right)  $ and
whose eigenvalues are $\left(  1,k,k+1,...,2k+1,3k\right)  .$ There exist
classifications of rigid algebras having similar sequences of eigenvalues and
filiform nilradical. However, there is nothing similar for
nonfiliform Lie algebras. Now the interesting fact is that we can extend centrally the preceding nilradicals of rigid laws to obtain characteristically nilpotent Lie algebras [22] : 

\begin{theorem}
The Lie algebras $e_{1}\left(  \frak{g}_{\left(  m,m-1\right)  }+\psi\right)
$ $\left(  m\geq3\right)  $ defined by the structural equations
\begin{align*}
d\omega_{1} &  =d\omega_{2}=0\\
d\omega_{j} &  =\omega_{1}\wedge\omega_{j-1},\;3\leq j\leq2m-1\\
d\omega_{2m} &  =\omega_{1}\wedge\omega_{2m-1}+\omega_{2}\wedge\omega_{3}\\
d\omega_{2m+1} &  =\sum_{j=2}^{\left[  \frac{2m+1}{2}\right]  }\left(
-1\right)  ^{j}\;\omega_{j}\wedge\omega_{2m+1-j}\\
d\omega_{2m+2} &  =\omega_{1}\wedge\omega_{2m+1}+\sum_{j=2}^{\left[
\frac{2m+1}{2}\right]  }\left(  -1\right)  ^{j}\,\left(  m+1-j\right)
\;\omega_{j}\wedge\omega_{2m+2-j}%
\end{align*}
are characteristically nilpotent.
\end{theorem}

\begin{corollary}
There are characteristically nilpotent Lie algebras $\frak{g}$ with nilindex
$2m+2$ for any $m\geq3.$
\end{corollary}

observe that the previous algebras have characteristic sequence $\left(  2m,1,1\right)  .$
This fact is directly related with the position of the vector $X_{2m+1}$ in
the graduation of $\frak{g}_{\left(  m,m-1\right) }$. The joined differential
form involves the form $\omega_{2}\wedge\omega_{2m}$, so that the nilindex of
the algebra increases. Moreover, observe that we have
\[
e_{1}\left(  \frak{g}_{\left(  m,m-1\right)  }+\psi\right)  \simeq
e_{1}\left(  \frak{g}_{\left(  m,m-1\right)  }\right)  +\psi
\]
so that we could have constructed the algebras extending and then deforming by
taking the same deformation. This gives, in a certain manner, a procedure to
generate characteristically nilpotent Lie algebras by extensions and
deformations of naturally graded Lie algebras ( see [10]).

\section{Affine structures over Lie algebras}

The origin of affine structures over Lie algebras is the tudy of affine
left-invariant structures over Lie groups [9]. The question wheter any
solvable Lie group admits a left invariant affine structure is a problem of
great interest, as it relates geometrical aspects of affine manifold theory
with representation theory of Lie algebras. Translated into Lie algebra
language, the question is if any solvable Lie algebra satisfies a certain
condition which is called affine structure. This goes back to Milnor in the
seventies, and is therefore called the Milnor conjeture. By the time the
problem was posted, all known results referred to low dimensions, where the
answer is positive. The first counterexample to Milnor' s conjecture was given
by Benoist [12]. He constructed a 11-dimensional filiform Lie algebra which
does not admit an affine structure. Explicitely, let $\frak{a}\left(  t\right)  $ be the filiform Lie algebra
given by
\begin{align*}
\left[  X_{1},X_{i}\right]    & =X_{i+1},\;2\leq i\leq10\\
\left[  X_{2},X_{3}\right]    & =X_{5}\\
\left[  X_{2},X_{5}\right]    & =-2X_{7}+X_{8}+tX_{9}%
\end{align*}
over the basis $\left\{  X_{1},..,X_{11}\right\}  $. The main point is to
prove that this algebra does not admit a faithful representation of degree 12,
which proves the nonexistence of an affine connection [12]. This example has
been widely generalized in [15] :  

\begin{theorem}
There exist filiform Lie algebras of dimensions $10\leq n\leq12$ which do not
admit an affine structure. For $n\leq9$ any filiform Lie algebra admits an
affine structure.
\end{theorem}

For this, cohomological methods are of importance, in particular the dimensions
of the cohomology spaces $H^{2}\left(  \frak{g},\mathbb{C}\right)  $, which
are usually called Betti numbers.
Let $\frak{g}$ be an $n$-dimensional Lie algebra and $G$ its associated Lie
group. If the group posseses a left-invariant affine structure, then this
induces a flat torsionfree left-invariant affine conection $\nabla$ on $G$,
that is
\begin{align*}
\nabla_{X}Y-\nabla_{Y}X-\left[ X%
,Y\right]    & =0\\
\nabla_{X}\nabla_{Y}Z-\nabla_{Y}\nabla_{X%
}Z-\nabla_{\left[ X,Y\right]  }Z  & =0
\end{align*}
for all left invariant vector fields $X,Y,Z$ on $G$. Now, defining%
\[
X.Y=\nabla_{X}Y%
\]
we obtain a bilinear product which satisfies
\[
X.\left(  Y.Z\right)  -\left(  X.Y\right)  .Z-Y.\left(  X.Z\right)  +\left(
Y.X\right)  .Z=0
\]
Observe that this implies that the product is left symmetric. 

\begin{definition}
An affine structure on a Lie algebra $\frak{g}$ is a bilinear product
$\frak{g}\times\frak{g}\rightarrow\mathbb{C}$ which is left symmetric and
satisfies
\[
\left[  X,Y\right]  =X.Y-Y.X,\;\forall X,Y\in\frak{g}%
\]
\end{definition}

It is known that there exists a one-to-one correspondence between affine
structures on $\frak{g}$ and left invariant structures on the associated Lie
group $G$ [28]. The interesting fact is that the problem can be dealt with
methods of representation theory of nilpotent Lie algebras.

\begin{proposition}
Let $\frak{g}$ be an $n$-dimensional Lie algebra. If $\frak{g}$ admits an
affine structure then $\frak{g}$ possesses a faithful module $M$ of dimension
$n+1$.
\end{proposition}

By the theorem of Ado [1], any Lie algebra admits a faithful representation.
Unfortunately, the results does not say anything about the minimal degree of
such a representation. Nowadays, it is accepted that the best lower bound is
given in [16]. This bound, equal to $\frac{\alpha}{\sqrt{n}}2^{n}$ with
$\alpha\sim2,76287$, has been used to obtain other counterexamples to Milnor's
conjecture [16]. In relation with the derivations structure, we have the
following 

\begin{proposition}
A Lie algebra $\frak{g}$ admits an affine structure if and only if there is a
$\frak{g}$-module $M$ of dimension $\dim\frak{g}$ such that $Z^{1}\left(
\frak{g},M\right)  $ contains a nonsingular cocycle.
\end{proposition}

The result is a consequence of the inversibility for a nonsingular cocycle. An
immediate corollary is 

\begin{corollary}
If $\frak{g}$ admits a nonsingular derivation, then it admits an affine structure.
\end{corollary}

Observe in particular the importance of this for graded Lie algebras : if
$\frak{g}$ is naturally graded ( the results remains valid for any positive
indexed graded Lie algebra ) then the natural operation defines a nonsingular
derivation, from which we obtain that any naturally graded Lie algebra has an
affine structure. As it is known that metabelian Lie algebras and those of
dimensions $n\leq6$ can be graded in such manner, all them admit an affine
structure. For $3$-step nilpotent Lie algebras Scheunemann [87] proved in
1974 the following :

\begin{theorem}
Any 3-step nilpotent Lie algebra $\frak{g}$ admits an affine structure. 
\end{theorem}

Observe that the algebra of Dixmier and Lister is 3-step nilpotent, thus it
has such a structure. Clearly all derivations are singular, which proves the
existence of CNLA with affine structures. The question is which of the
structural properties of CNLA allow the existence of such structures. In
particular, has it any relation with the structure of the automorphism group
?\newline For 4-step nilpotent Lie algebras the question is open, and the bEst
result achieved can be found in [27]. However, the fundamental source ( once
more ) for the study of affine structures is the variety $\mathcal{F}^{n}$ of
filiform Lie algebras. In [16] the author defines the following cocycles :

\begin{definition}
Let $\frak{g}$ be a filiform Lie algebra. A cocycle $\omega\in Z^{2}\left(
\frak{g},\mathbb{C}\right)  $ is called affine if it is nonzero over $Z\left(
\frak{g}\right)  \wedge\frak{g}$. A class $\left[  \omega\right]  \in
H^{2}\left(  \frak{g},\mathbb{C}\right)  $ is called affine if every
representative is affine.
\end{definition}

Then the next result characterizes certain extensions of filfirom Lie algebras :

\begin{proposition}
A filiform Lie algebra $\frak{g}$ has a filiform extension of dimension
dim$\frak{g}$+1 if and only if there exists an affine cohomology class in
$\frak{g}$. 
\end{proposition}

This result has two interesting consequences :

\begin{corollary}
If the filiform Lie algebra $\frak{g}$ admits an affine cohomology class
$\left[  \omega\right]  $, then it admits an affine structure.
\end{corollary}

\begin{corollary}
If $\frak{g}$ is filiform of dimension $n\geq6$ and $\dim H^{2}\left(
\frak{g},\mathbb{C}\right)  =2$, then $\frak{g}$ has no affine cohomolofgy class.
\end{corollary}

Endowed with these methods, Burde has constructed two classes of filiform Lie
algebras [17] which provide a lot of counterexamples to Milnors conjecture.

We conclude giving CNLAs which admit an affine structure but whose Lie algebra
of derivations is not characteristically nilpotent : over the basis $\left\{
X_{1},..,X_{11}\right\}  $ let $\frak{g}_{\left(  a_{i}\right)  }$ be the
filiform Lie algebra given by
\[%
\begin{array}
[c]{ll}%
\left[  X_{1},X_{i}\right]  =X_{i+1},\;2\leq i\leq10 & \\
\left[  X_{2},X_{3}\right]  =X_{5} & \left[  X_{3},X_{6}\right]  =-\frac
{12}{5}X_{9}+a_{5}X_{10}+a_{6}X_{11}\\
\left[  X_{2},X_{4}\right]  =X_{6} & \left[  X_{3},X_{7}\right]  =-\frac
{39}{5}X_{10}+a_{7}X_{11}\\
\left[  X_{2},X_{5}\right]  =-2X_{7}+X_{8} & \left[  X_{3},X_{8}\right]
=a_{8}X_{11}\\
\left[  X_{2},X_{6}\right]  =-5X_{8}+2X_{9} & \left[  X_{4},X_{5}\right]
=\frac{27}{5}X_{9}+a_{9}X_{10}+a_{10}X_{11}\\
\left[  X_{2},X_{7}\right]  =-\frac{13}{9}X_{9}+a_{1}X_{10}+a_{2}X_{11} &
\left[  X_{4},X_{6}\right]  =\frac{27}{5}X_{10}+a_{9}X_{10}\\
\left[  X_{2},X_{8}\right]  =\frac{26}{5}X_{10}+a_{3}X_{11} & \left[
X_{4},X_{7}\right]  =a_{11}X_{11}\\
\left[  X_{2},X_{9}\right]  =a_{4}X_{11} & \left[  X_{5},X_{6}\right]
=a_{12}X_{11}\\
\left[  X_{3},X_{4}\right]  =3X_{7}-X_{8} & \\
\left[  X_{3},X_{5}\right]  =3X_{8}-X_{9} &
\end{array}
\]
where the following relations are satisfied :%
\[
\left\{
\begin{array}
[c]{cc}%
a_{5}-3a_{1}+\frac{26}{5}-a_{9}=0; & 3a_{2}-a_{3}+a_{10}-a_{6}=0\\
a_{3}-3a_{2}-2a_{10}=0; & a_{12}-\frac{27}{5}a_{4}-\frac{54}{5}=0\\
2a_{2}-a_{3}+2a_{6}-a_{5}-a_{7}+2=0; & -3a_{3}+a_{4}+2+4a_{9}-\frac{31}%
{5}-a_{11}=0\\
a_{1}+a_{5}-2=0; & 5a_{3}-2a_{4}-2a_{8}+\frac{52}{2}+5a_{5}+5a_{7}-10=0\\
a_{7}+a_{3}-a_{1}=0; & a_{4}+a_{8}-\frac{26}{5}=0\\
a_{9}-4a_{5}+6-\frac{26}{5}=0 &
\end{array}
\right\}
\]

This example, for the values
\begin{align*}
a_{1}  & =\frac{51}{25};\;a_{2}=-a_{6}=a_{10}=\frac{28}{125};\;a_{3}=\frac
{28}{25};\;a_{4}=\frac{19}{16};\;a_{5}=-\frac{1}{25}\\
a_{7}  & =\frac{23}{25};\;a_{8}=\frac{321}{80};\;a_{9}=-\frac{24}{25}%
;\;a_{11}=-\frac{189}{16};\;a_{12}=\frac{1377}{80}%
\end{align*}
is due to Remm and Goze.   

\section{Associative characteristic nilpotent algebras}

Motivated by the paper of Dixmier and Lister, in 1971 T.\ S. Ravisankar [81]
extended the concept of being characteristically nilpotent to general
algebras. \ This approach has been useful for the study of Malcev algebras, as
for associative algebras and its deformation theory [72]. 

Let $A$ be a nonassociative complex algebra ( again we convene that the base
field is $\mathbb{C}$, though this assumption is not generally necessary ). We
denote its Lie algebra of derivations by $D\left(  A\right)  $. Let
\[
A^{\left[  1\right]  }=\left\{  \sum D_{i}x_{i}\;|\;x_{i}\in A,\;D_{i}\in
D\left(  A\right)  \right\}
\]
and define inductively $A^{\left[  k+1\right]  }=$ $\left\{  \sum D_{i}%
y_{i}\;|\;y_{i}\in A^{[k]},\;D_{i}\in D\left(  A\right)  \right\}  $. 

\begin{definition}
An algebra $A$ is called characteristically nilpotent ( C-nilpotent ) if there
exists an integer $n$ such that $A^{[n]}=0$.
\end{definition}

It is clear that if $A$ is a C-nilpotent algebra, then any derivation of $A$
is a linear nilpotent transformation on $A$. The converse also holds [81].
For the special case of associative algebras, in [43], let $e_{\alpha}$ be
the $\left(  r+1\right)  ^{2}$ matrix whose $\alpha=\left(  i,j\right)  $
entry is one, otherwise zero. The space generated by this vector is denoted by
$E_{\alpha}$. Let us then define%
\begin{align*}
R  & =\left\{  \alpha=\left(  i,j\right)  ,\;1\leq i,j\leq r+1\right\}  \\
R^{+}  & =\left\{  \alpha=\left(  i,j\right)  \;|\;i<j\right\}
\end{align*}
It follows that $R^{-}=R-R^{+}$ and $S=\left\{  \left(  i,i+1\right)  \in
R\;|\;1\leq i\leq r\right\}  $ is the set of simple roots, in analogy with the
Lie algebra case [23]. Then  $L=\sum_{R^{+}}E_{\alpha}$ is a nilpotent
associative algebra. Consider the bilinear mappings of $L\times L\rightarrow
L$ defined by
\[
g_{k,m}\left(  e_{\alpha_{k}},e_{\alpha_{m}}\right)  =e_{\delta},\;\text{where
}\alpha_{i}=\left(  i,i+1\right)  \text{ and }\delta=\left(  1,r+1\right)
\]
Obviously the center of $L$ is generated by the root $\delta$. Let us now
consider the linear combination $\psi=\sum_{1\leq k,m\leq r}a_{k,m}g_{k,m}$
$\ $for $a_{k,m}\in\mathbb{C}$. In [43] it is proven that this is a lineraly
expandable cocycle, and further that

\begin{theorem}
Let $\psi$ be the cocycle given by $\psi=\sum_{1\leq k,m\leq r}a_{k,m}g_{k,m}$
$\ $with $\prod_{1\leq i\leq r}a_{ii}\neq0$. Then the associative algebra
$L+\psi$ is characteristically nilpotent.
\end{theorem}

Constructing families of this kind, the variety $\mathcal{N}^{n}$ of
associative algebras can be studied as Lie algebras have been [72]. In
particular, among other results the following shows the similarity between the
theory of characteristically nilpotent Lie algebras and C-algebras :

\begin{theorem}
For $n\geq2$ there exists a Zariski-open subset of $\mathcal{N}^{n}$ formed by
characteristically nilpotent associative algebras. Moreover, its dimension is
$n^{2}-n$.
\end{theorem}

\end{document}